\documentclass[11pt]{amsart}
\usepackage{amsmath,amsthm,amsfonts,amscd,amssymb,eucal,latexsym,mathrsfs}
\usepackage[all]{xy}

\theoremstyle{definition}
\newtheorem{theorem}{Theorem}[section]
\newtheorem{corollary}[theorem]{Corollary}
\newtheorem{lemma}[theorem]{Lemma}
\newtheorem{proposition}[theorem]{Proposition}
\newtheorem{remark}[theorem]{Remark}

\newtheorem{example}[theorem]{Example}
\newtheorem{question}[theorem]{Question}
\newtheorem{definition}[theorem]{Definition}

\newcommand{\htopol}{{\text{\rm h}}_{\text{\rm top}}}
\newcommand{\rank}{{\rm rank}}

\newcommand{\sa}{{\rm sa}}
\newcommand{\Aff}{{\rm Aff}}
\newcommand{\sep}{{\rm sep}}
\newcommand{\spa}{{\rm spn}}
\newcommand{\spn}{{\rm span}}

\newcommand{\Aut}{{\rm Aut}}
\newcommand{\interior}{{\rm int}}
\newcommand{\CA}{{\rm CA}}
\newcommand{\rc}{{\rm rc}}
\newcommand{\Fin}{\mathcal{P}_{\rm f}}
\newcommand{\hc}{{\rm hc}}
\newcommand{\hcc}{{\rm hcc}}
\newcommand{\Ad}{{\rm Ad}\,}

\newcommand{\id}{{\rm id}}

\newcommand{\Inn}{{\rm Inn}}
\newcommand{\IA}{{\rm IA}}
\newcommand{\entr}{{\rm ht}}

\newcommand{\CCA}{{\rm CCA}}
\newcommand{\rcc}{{\rm rcc}}

\newcommand{\Prim}{{\rm Prim}}

\newcommand{\VB}{{\rm ht}}
\newcommand{\red}{{\rm r}}

\newcommand{\Asp}{{\rm Asp}}

\begin{document}

\title[Dynamical entropy in Banach spaces]{Dynamical entropy in Banach spaces}

\author{David Kerr}
\author{Hanfeng Li}
\address{\hskip-\parindent
David Kerr, Department of Mathematics, Texas A{\&}M University, 
College Station TX 77843-3368, U.S.A. and
Graduate School of Mathematical Sciences, University 
of Tokyo, 3-8-1 Komaba, Meguro-ku, Tokyo 153-8914, Japan}
\email{kerr@ms.u-tokyo.ac.jp}

\address{\hskip-\parindent
Hanfeng Li, Department of Mathematics, University of Toronto, Toronto, Ontario M5S 3G3, 
Canada}
\email{hli@fields.utoronto.ca}

\date{April 12, 2005}

\begin{abstract}
We introduce a version of Voiculescu-Brown approximation entropy for isometric
automorphisms of Banach spaces and develop within this framework the 
connection between dynamics and the local theory of Banach spaces as 
discovered by Glasner and Weiss. Our fundamental result concerning this 
contractive approximation entropy, or CA entropy, characterizes the 
occurrence of positive values both geometrically and topologically. This 
leads to various applications; for example, we obtain a geometric description
of the topological Pinsker factor and show that a $C^*$-algebra is type I
if and only if every multiplier inner $^*$-automorphism has zero CA entropy.
We also examine the behaviour of CA entropy under various product
constructions and determine its value in many examples, including isometric
automorphisms of $\ell_p$ for $1\leq p \leq\infty$ and noncommutative tensor
product shifts.
\end{abstract}

\maketitle

\section{Introduction}

In \cite{GW} E. Glasner and B. Weiss showed that if
a homeomorphism from a compact metric space $K$ to itself has zero
topological entropy, then so does the induced homeomorphism on the space
of probability measures on $K$ with the weak$^*$ topology. One of the
two proofs they gave of this striking result established a remarkable
connection between topological dynamics and the local
theory of Banach spaces. The key geometric fact is the exponential
dependence of $k$ on $n$ given an
approximately isometric embedding of $\ell^n_1$ into $\ell^k_\infty$,
which is a consequence of the work of T. Figiel, J. Lindenstrauss, and
V. D. Milman on almost Hilbertian sections of unit balls in Banach
spaces \cite{FLM}.

The first author showed in \cite{EID} that Glasner and Weiss's
geometric approach can be conceptually simplified from a
functional-analytic viewpoint using Voiculescu-Brown entropy and also
more generally applied to show that if a $^*$-automorphism of a separable
exact $C^*$-algebra has zero Voiculescu-Brown entropy then the induced
homeomorphism on the unit ball of the dual has zero topological entropy.
In this case the crucial Banach space fact is the exponential dependence of
$k$ on $n$ given an approximately isometric embedding of $\ell^n_1$ into
the matrix $C^*$-algebra $M_k$ \cite[Lemma 3.1]{EID}, which can be deduced
from the work of N. Tomczak-Jaegermann on the Rademacher type $2$ constants
of Schatten $p$-classes \cite{TC}.

In the present paper we pursue this connection between dynamics and Banach
space geometry within a general Banach space framework via the introduction
of an analogue of Voiculescu-Brown entropy which we call
contractive approximation entropy, or simply CA entropy. This dynamical
invariant has the advantage of being defined for any isometric automorphism of
a Banach space and often exhibits greater tractability as a result of its more
basic structural context. With entropy a basic problem is to
determine whether or not it is positive (in which case the dynamics can be
thought of as being chaotic, or nondeterministic), and in the case of CA
entropy we are able to characterize the occurrence of positive values both
geometrically and topologically by expanding upon the arguments of
Glasner and Weiss. A large part of the paper
will involve applications of these characterizations. In particular, we give
a description of the Pinsker algebra in topological dynamics in terms of
dynamically generated sets equivalent to the standard basis of $\ell_1$, and
prove that a $C^*$-algebra is type I if and only if every multiplier inner
$^*$-automorphism has zero CA entropy. The latter result was conjectured by
N. P. Brown for Voiculescu-Brown entropy \cite{typeI} but seems to be out of
reach in that case. For $C^*$-dynamics, the drawback of CA entropy in
comparison to Voiculescu-Brown entropy is that, being oblivious to 
matricial structure, it can be much cruder as a numerical invariant, as
we illustrate by computing the value for the tensor product shift on the
UHF algebra $M_p^{\otimes\mathbb{Z}}$
to be infinity. Nevertheless, the condition of
simplicity on a $C^*$-algebra does not necessarily rule out the existence of
$^*$-automorphisms with finite nonzero CA entropy, as we have discovered with
the type $2^\infty$ Bunce-Deddens algebra.

We begin in Section~\ref{S-CA} by defining CA entropy and recording some basic
properties. In Section~\ref{S-zero}
we demonstrate that the topological entropy of a homeomorphism of compact
Hausdorff space agrees with the CA entropy of the induced $^*$-automorphism of
the $C^*$-algebra of functions over the space, and then proceed to establish
our geometric and topological characterizations of positive CA entropy
(Theorem~\ref{T-zero}).
Specifically, we show that positive CA entropy is equivalent to positivity of
the topological entropy of the induced homeomorphism on the unit ball of the
dual as well as to the existence of an element the restriction of whose
orbit to a positive density subset of iterates is equivalent to the standard
basis of $\ell_1$. Several immediate corollaries ensue,
such as the invariance of positive CA entropy under isomorphic conjugacy and
the vanishing of CA entropy for every isometric automorphism of a Banach space
with separable dual. We round out Section~\ref{S-zero} by exhibiting
$C(\mathbb{T} )$ as an example of a Banach space which contains $\ell_1$
isometrically but admits no isometric automorphism with positive CA entropy.

Section~\ref{S-mat} briefly compares CA entropy with its matricial analogue
for completely isometric automorphisms of exact operator spaces,
which we call completely contractive approximation entropy, or CCA entropy.
This was introduced in the $C^*$-algebraic setting by C. Pop and R. Smith, who
showed that it coincides with Voiculescu-Brown entropy in that case \cite{PS}.
Positive CA entropy implies positive CCA entropy, but the
converse is false. However, we have been unable to resolve the problem of the
converse for $^*$-automorphisms of exact $C^*$-algebras. As our computation
for the noncommutative shift in Section~\ref{S-shift} demonstrates, it is also
possible for the CA entropy to be strictly larger than the CCA entropy.

In Section~\ref{S-tPa} we obtain a geometric description of the topological
Pinsker algebra, which is the $C^*$-algebraic manifestation of the
Pinsker factor (i.e., the largest zero entropy factor) in topological
dynamics \cite{ZEF}. It turns out that the elements of the topological
Pinsker algebra are precisely those which do not generate a subspace
canonically isomorphic to $\ell_1$ along a positive density set of
iterates. This yields geometric characterizations of positive and
completely positive topological entropy, and implies that tame
and HNS $\mathbb{Z}$-systems \cite{tame,GM} have zero entropy.

Section~\ref{S-typeI} contains our entropic characterization of type I
$C^*$-algebras. For this we establish two lemmas which demonstrate that
the property of having zero CA entropy behaves well with respect to
Banach space quotients and continuous fields of Banach spaces.

Section~\ref{S-prod} examines the behaviour of CA entropy under various
product constructions. Subadditivity holds for the injective tensor product,
but not necessarily for other tensor product norms, as we illustrate with the
shift on a spin system. Adapting an argument of C. Pop and
R. Smith based on Imai-Takai duality \cite{PS}, we prove that, under taking a
$C^*$-algebraic crossed product by the action of an amenable locally compact
group $G$, zero CA entropy is preserved on a group element $g$ such that
$\Ad g$ generates a compact subgroup of $\Aut (G)$, and if the original
$C^*$-algebra is commutative then every value of CA entropy is preserved
in the same situation. The latter fact has the consequence that there exist
simple $C^*$-algebras, for instance the type $2^\infty$ Bunce-Deddens algebra,
which admit inner $^*$-automorphisms with any prescribed value of CA entropy.
Finally we deduce in Section~\ref{S-prod} that zero entropy is preserved under
taking reduced free products of commutative probability spaces by applying a 
result of N. P. Brown, K. Dykema, and D. Shlyakhtenko \cite{BDS}.

In Section~\ref{S-opalg} we study the prevalence of zero and infinite
CA entropy in $C^*$-algebras. We find that, for many $C^*$-algebras which
are subject to classification theory (more specifically, those which are
tensorially stable with respect to the Jiang-Su algebra $\mathcal{Z}$), the
collection of $^*$-automorphisms
with infinite CA entropy is point-norm dense, while in the special cases of UHF
algebras, the Cuntz algebra $\mathcal{O}_2$, and the Jiang-Su algebra
the collection of $^*$-automorphisms with zero CA entropy is a point-norm
dense $G_\delta$ set, giving a noncommutative version of a result of Glasner
and Weiss on homeomorphisms of the Cantor set \cite{TRP}.

In the last three sections we apply a combinatorial argument 
(Lemma~\ref{L-comb}) to compute the CA entropy for some canonical examples.
In Section~\ref{S-shift} we prove that the tensor product shift on the UHF 
algebra $M_p^{\otimes\mathbb{Z}}$ has infinite CA entropy, in contrast to the 
value of $\log p$ for the Voiculescu-Brown entropy. This raises the question
of whether there exist $^*$-automorphisms of simple AF algebras with finite
nonzero CA entropy. In Section~\ref{S-ellinfty} we show
that the CA entropy
of an isometric automorphism of $\ell_\infty$ is either zero or infinity
depending on whether or not there is a finite bound on the cardinality of the
orbits of the associated permutation of $\mathbb{Z}$.
Finally, in Section~\ref{S-ell1} we show that the CA entropy
of an isometric automorphism of $\ell_1$ is either zero or infinity
depending on whether or not there is an infinite orbit in the
associated permutation of $\mathbb{Z}$. As a consequence we deduce that CA 
entropy is not an isomorphic conjugacy invariant.

All Banach spaces will be over the complex numbers, unless there is an
indication to the contrary, such as the tag $\mathbb{R}$ when referring to
the real version of a standard Banach space. We point out however that the
relevant results in Sections~\ref{S-CA}, \ref{S-zero}, and
\ref{S-ellinfty} are also valid over
the real numbers, although in other situations differences between the
real and complex cases can arise, as Remark~\ref{R-sphere} illustrates.
For terminology related to Banach spaces see \cite{CBS,JL}. The spaces
$\ell_p$ for $1\leq p\leq\infty$ will be indexed over $\mathbb{Z}$ so that
we may conveniently speak of the shift automorphism as obtained from the
shift $k\mapsto k+1$ on $\mathbb{Z}$. Given a Banach space $X$
and $r>0$ the ball $\{ x\in X : \| x \| \leq r \}$ will be denoted $B_r (X)$.
All $C^*$-tensor products will be minimal and written using an
unadorned $\otimes$. The set of self-adjoint elements of
an operator system $X$ will be denoted $X_\sa$. We write $M_d$ for the
$C^*$-algebra of $d\times d$ complex matrices.

This paper subsumes the material from our e-print \cite{PTE}.
\medskip

\noindent{\it Acknowledgements.} D. Kerr was supported by the
Alexander von Humboldt Foundation and heartily thanks Joachim Cuntz for
hosting his stay at the University of M\"{u}nster over the 2003--2004 academic
year. The preliminary stages of this work were carried out during his
stay at the University of Rome ``La Sapienza'' over the 2002--2003 academic
year and he is grateful to Claudia Pinzari for her generous
hospitality and to NSERC for support. He thanks Vitali Milman for
addressing a question in the local theory of Banach spaces. H. Li thanks
George Elliott for helpful discussions. We thank Andrew
Toms and Wilhelm Winter for discussions relating to $\mathcal{Z}$-stability.
We are also indebted to the referee for suggesting a unified 
approach to obtaining the lower bounds in 
Sections~\ref{S-shift} to \ref{S-ell1} via
Lemma~\ref{L-comb}, a variant of one of our original lemmas.

\section{Contractive approximation entropy}\label{S-CA}

Let $X$ and $Y$ be Banach spaces and $\gamma : X\to Y$ a bounded linear map.
Denote by $\Fin (X)$ the collection of finite subsets of $X$.
For each $\Omega\in\Fin (X)$ and $\delta > 0$ we denote by
$\CA (\gamma , \Omega , \delta )$ the collection of triples
$(\phi , \psi ,d)$ where
$d$ is a positive integer and $\phi : X\to\ell^d_\infty$ and
$\psi : \ell^d_\infty \to Y$ are contractive linear maps such that
$$ \| \psi\circ\phi (x) - \gamma (x) \| < \delta $$
for all $x\in\Omega$. By a {\it CA embedding} of a Banach space $X$ we mean
an isometric linear map $\iota$ from $X$ to a Banach space $Y$ such that
$\CA (\iota , \Omega , \delta )$
is nonempty for every $\Omega\in\Fin (X)$ and $\delta > 0$.
Every Banach space admits a CA embedding; for example, the canonical
map $X\to C(B_1 (X^* ))$ defined via evaluation is a
CA embedding, as a standard partition of unity argument shows.

Let $\iota : X\to Y$ be a CA embedding. For each $\Omega\in\Fin (X)$
and $\delta > 0$ we set
$$ \rc (\Omega , \delta ) = \inf \{ d : (\phi , \psi ,d)\in
\CA (\iota , \Omega , \delta ) \} . $$
We claim that this quantity is independent of the CA embedding, as our
notation indicates. Indeed
suppose $\iota_0 : X\to Y_0$ is another CA embedding and
$(\phi , \psi ,d)\in \CA (\iota , \Omega , \delta )$. Take an
$\varepsilon > 0$ such that
$$ \| \psi\circ\phi (x) - \iota (x) \| < \delta - \varepsilon $$
for all $x\in\Omega$, and choose a $(\phi_0 , \psi_0 ,d_0 )\in
\CA (\iota_0 , \Omega , \varepsilon )$.
By the injectivity of $\ell^{d_0}_\infty$ we can extend $\phi\circ\iota^{-1}
: \iota (X) \to\ell^{d_0}_\infty$ to a contractive linear map
$\rho : Y\to\ell^{d_0}_\infty$. An application of the triangle inequality
then shows that $(\phi , \psi_0 \circ\rho\circ\psi , d)\in
\CA (\iota_0 , \Omega , \delta )$,
from which the claim follows.

We denote by $\IA (X)$ the collection of all isometric automorphisms of $X$.
For $\alpha\in\IA (X)$ we set
\begin{align*}
\hc (\alpha , \Omega , \delta ) &= \limsup_{n\to\infty} \frac1n
\log\rc (\Omega\cup\alpha\Omega \cup\cdots\cup
\alpha^{n-1} \Omega , \delta ) , \\
\hc (\alpha , \Omega ) &= \sup_{\delta > 0}
\hc (\alpha , \Omega , \delta ) , \\
\hc (\alpha ) &= \sup_{\Omega\in\Fin (X)} \hc (\alpha , \Omega )
\end{align*}
and refer to the last quantity as the {\it contractive approximation
entropy} or {\it CA entropy} of $\alpha$.

From the fact that $\rc (\Omega , \delta )$ does not depend
on the CA embedding we see that CA entropy is an invariant with respect to
conjugacy by isometric isomorphisms. It is not, however, an invariant with
respect to conjugacy by arbitrary isomorphisms (see Remark~\ref{R-conj}).
Nevertheless, it turns out that zero CA entropy is an isomorphic conjugacy
invariant, as we will see in Section~\ref{S-zero}.

\begin{proposition}\label{P-prop}
Let $\alpha$ be an isometric automorphism of a Banach space $X$.
\begin{enumerate}
\item[(i)] If $Y\subseteq X$ is an $\alpha$-invariant closed
subspace then $\hc (\alpha |_Y ) \leq \hc (\alpha )$ (monotonicity).

\item[(ii)] If $\{ \Omega_\lambda \}_{\lambda\in\Lambda}$ is an increasing net
in $\Fin (X)$ such that $\bigcup_{\lambda\in\Lambda} \bigcup_{n\in\mathbb{Z}}
\alpha^n (\Omega_\lambda )$ is total in $X$ then
$\hc (\alpha ) = \sup_\lambda \hc (\alpha , \Omega_\lambda )$,

\item[(iii)] For every $k\in\mathbb{Z}$ we have
$\hc (\alpha^k ) = |k| \hc (\alpha )$.
\end{enumerate}
\end{proposition}

\begin{proof}
Monotonicity follows from the fact that the restriction of a CA embedding to a
closed subspace is a CA embedding, while for (ii) and (iii) we can proceed as
in the proofs of Proposition~1.3 and 3.4, respectively, in \cite{Voi}.
\end{proof}

\begin{proposition}\label{P-directsum}
Let $X_1 , \dots , X_r$ be Banach spaces with respective isometric
automorphisms $\alpha_1 , \dots , \alpha_r$. Then for the isometric
automorphism $\alpha_1 \oplus\cdots\oplus \alpha_r$ of the $\ell_\infty$-direct
sum $(X_1 \oplus\cdots\oplus X_r )_\infty$ we have
$$ \hc (\alpha_1 \oplus\cdots\oplus \alpha_r ) = \max_{1\leq i \leq r}
\hc (\alpha_i ) . $$
\end{proposition}

\begin{proof}
The inequality $\hc (\alpha_1 \oplus\cdots\oplus \alpha_r ) \geq
\max_{1\leq i \leq r} \hc (\alpha_i )$ is a consequence of monotonicity
(Proposition~\ref{P-prop}(i)), while the reverse inequality is readily
seen using the fact that an
$\ell_\infty$-direct sum of CA embeddings is a CA embedding.
\end{proof}

\noindent We also see by applying Proposition~\ref{P-directsum} in conjunction
with Proposition~\ref{P-prop}(ii) that the CA entropy of a $c_0$-direct sum of
isometric automorphisms is equal to the supremum of the CA entropies of the
summands.

\section{Topological and geometric characterizations of
positive CA entropy}\label{S-zero}

Let $K$ be a compact Hausdorff space and $T : K\to K$ a homeomorphism. Recall
that the {\it topological entropy} of $T$, denoted $\htopol (T)$, is defined
as the supremum over all finite open covers $\mathcal{W}$ of $K$ of the
quantities
$$ \lim_{n\to\infty} \frac1n \log
N(\mathcal{W}\vee T^{-1}\mathcal{W}\vee\cdots\vee T^{-(n-1)} \mathcal{W}) , $$
where $N(\cdot )$ denotes the smallest cardinality of a subcover \cite{AKM}
(see \cite{DGS,HK,Wal} for general references). The topological entropy of
$T$ may be equivalently expressed in terms of separated and spanning
sets \cite{Bow,Mis} as follows.

Denote by $\mathcal{U}$ the unique uniformity compatible with the
topology on $K$, i.e., the collection of all neighbourhoods of the diagonal
in $K\times K$. Let $Q\subseteq K$ be a compact subset, and let
$U\in\mathcal{U}$. A set $E\subseteq K$ is
{\it $(n,U)$-separated (with respect to $T$)} if for every $s,t\in
E$ with $s\neq t$ there exists a $0\leq k \leq n-1$ such that $(T^k s , T^k t)
\notin U$. A set $E\subseteq K$ is {\it $(n,U)$-spanning for
$Q$ (with respect to $T$)}
if for every $s\in Q$ there is a $t\in E$ such that $(T^k s , T^k t) \in U$
for each $k=0, \dots , n-1$. Denote by
$\sep_n (T,Q,U)$ the largest cardinality of an
$(n,U)$-separated subset of $Q$ and
by $\spa_n (T,Q,U)$ the smallest cardinality of an
$(n,U)$-spanning set for $Q$. When $Q=K$ we simply write
$\sep_n (T,U)$ and $\spa_n (T,U)$. We then have, by the same argument as
in \cite{Bow},
$$ \sup_{U\in\mathcal{U}} \limsup_{n\to\infty}\frac1n \log\sep_n (T,Q,U) =
\sup_{U\in\mathcal{U}} \limsup_{n\to\infty}\frac1n \log\spa_n (T,Q,U) . $$
We also obtain the same quantity by taking either of the suprema over any
given base for $\mathcal{U}$, the most important of which for our purposes
is the collection consisting of the sets
$U_{d,\varepsilon} = \{ (s,t) \in K\times K : d(s,t) < \varepsilon \}$ where
$d$ is a continuous pseudometric on $K$ and $\varepsilon > 0$.
This quantity, which is invariant under conjugacy by homeomorphisms, we denote
by $\htopol (T,Q)$. When $Q=K$ it can be shown, as in \cite{Bow}, that
we recover the topological entropy $\htopol (T)$.

In the following proposition $\VB (\cdot )$ denotes Voiculescu-Brown entropy
\cite{Voi,Br}.

\begin{proposition}\label{P-topol}
Let $K$ be a compact Hausdorff space and $T: K\to K$ a homeomorphism. Let
$\alpha_T$ be the $^*$-automorphism of $C(K)$ given by $\alpha_T (f) = 
f\circ T$ for all $f\in C(K)$. Then $\VB (\alpha_T ) = \hc (\alpha_T ) = 
\htopol (T)$.
\end{proposition}

\begin{proof}
The identity map from $C(K)$ to itself is a CA embedding, and the
partition of unity argument in the proof of Proposition~4.8 in \cite{Voi}
demonstrates that $\hc (\alpha_T ) \leq \htopol (T)$. The inequality
$\VB (\alpha_T ) \leq \hc (\alpha_T )$ follows from the formulation of
Voiculescu-Brown entropy in terms of completely contractive linear maps
\cite{PS} and the observation that the identity on $C(K)$ is a nuclear
embedding (see Section~\ref{S-mat}) along with the fact that a contractive
linear map from an operator space into a commutative $C^*$-algebra is
automatically completely contractive \cite[Thm.\ 3.8]{Pau}. It thus remains
to show that $\htopol (T) \leq \VB (\alpha_T )$. This holds when $K$ is
metrizable by Proposition~4.8 of \cite{Voi}, and we can reduce the
general case to the metrizable one as follows.

Let $\delta > 0$. Then there is a neighbourhood $U$ of the diagonal in
$K\times K$ such that
$$ \htopol (T) \leq \limsup_{n\to\infty}\frac1n \log\sep_n (T,U) + \delta . $$
We may assume that $U$ is of the form $U_{d,\varepsilon}$ for some
$\varepsilon > 0$ and pseudometric $d$ of the form
$$ d(s,t) = \sup_{f\in\Omega} | f(s) - f(t) | $$
for some $\Omega\in\Fin (C(K))$. Let $A$ be the unital $C^*$-subalgebra of
$C(K)$ generated by $\bigcup_{k\in\mathbb{Z}} \alpha^k (\Omega )$. Then $A$
is separable and $\alpha$-invariant. By separability we can enlarge $\Omega$
to a compact total subset $\Gamma$ of the unit ball of $A$ and define a
metric
$$ d' (s,t) = \sup_{f\in\Gamma} | f(s) - f(t) | $$
on the spectrum of $A$ which is compatible with the weak$^*$ topology, and with
respect to this metric the homeomorphism $S$ of the spectrum of $A$ induced
from $\alpha |_A$ evidently satisfies
$\sep_n (S,U_{d' ,\varepsilon} ) \geq \sep_n (T,U_{d,\varepsilon} )$
for all $n\in\mathbb{N}$. Knowing that the desired inequality holds
in the metrizable case and applying monotonicity we thus obtain
$$ \htopol (T) \leq \htopol (S) + \delta \leq \VB (\alpha |_A ) + \delta \leq
\VB (\alpha ) + \delta , $$
completing the proof.
\end{proof}

\begin{lemma}\label{L-ell1}
Let $X$ be a Banach space. Let $\Omega = \{ x_1 , \dots , x_n \} \subseteq X$ 
and suppose that the
linear map $\gamma : \ell_1^n \to X$ sending the $i$th standard basis element
of $\ell_1^n$ to $x_i$ for each $i=1, \dots ,n$ is an isomorphism.
Let $\delta > 0$ be such that $\delta < \| \gamma^{-1} \|^{-1}$. Then
$$ \log\rc (\Omega , \delta ) \geq na \| \gamma \|^{-2}
(\| \gamma^{-1} \|^{-1} - \delta )^2 $$
where $a>0$ is a universal constant.
\end{lemma}

\begin{proof}
Let $\iota : X\to Y$ be a CA embedding, and suppose $(\phi , \psi , d)
\in\CA (\iota , \Omega , \delta )$. For any linear combination
$\sum c_i x_i$ of the elements $x_1 , \dots , x_n$ we have
\begin{align*}
\Big\| \sum c_i x_i \Big\| &\leq \Big\| \iota\Big( \sum c_i x_i \Big)
- (\psi\circ\phi )\Big( \sum c_i x_i \Big) \Big\| \\
&\hspace*{35mm}\ + \Big\| (\psi\circ\phi )\Big( \sum c_i x_i \Big) \Big\| \\
&\leq \delta \sum |c_i | + \Big\| \phi \Big( \sum c_i x_i \Big) \Big\| \\
&\leq \delta \| \gamma^{-1} \|
\Big\| \sum c_i x_i \Big\| + \Big\| \phi \Big( \sum c_i x_i \Big) \Big\|
\end{align*}
and so $\big\| \phi \big( \sum c_i x_i \big) \big\| \geq (1 - \delta
\| \gamma^{-1} \| ) \big\| \sum c_i x_i \big\|$. Since $\phi$ is contractive,
it follows that the composition $\phi\circ\gamma$ is a
$\| \gamma \| (\| \gamma^{-1} \|^{-1} - \delta )^{-1}$-isomorphism
onto its image in $\ell_\infty^d$. The desired conclusion now follows
from the $\ell_\infty^d$ version of Lemma~3.1 in \cite{EID}, which can be
deduced by the same kind of argument using (Rademacher) type $2$ constants,
or as an immediate consequence by viewing $\ell_\infty^d$ as the diagonal in
the matrix algebra $M_d$.
\end{proof}

\noindent We remark in passing that Lemma~\ref{L-ell1} shows
that the CA entropy is infinite for the universal separable unital
$C^*$-dynamical system, i.e., the shift on the infinite full free product
$(C(\mathbb{T} )^{*\mathbb{N}} )^{*\mathbb{Z}}$, since the set of
canonical unitary generators is isometrically isomorphic to the
standard basis of $\ell_1$ (cf.\ \cite[Sect.\ 8]{Pis}).

The following lemma is well known in Banach space theory and follows readily
from Th\'{e}or\`{e}me 5 of \cite{Pajor}, as indicated in the proof of
Lemma~3.2 in \cite{EID}.

\begin{lemma}\label{L-proj}
For every $\varepsilon > 0$ and $\lambda > 0$ there exist $d > 0$ and
$\delta > 0$ such that the following holds for all $n\in\mathbb{N}$: if
$S\subseteq B_1 (\ell^n_{\infty} )$ is a symmetric convex set which contains
an $\varepsilon$-separated set $F$ of cardinality at least $e^{\lambda n}$
then there is a subset $I_n \subseteq \{ 1, 2, \cdots, n \}$ with cardinality
at least $dn$ such that
$$ B_{\delta} (\ell^{I_n}_{\infty} )\subseteq \pi_n (S), $$
where $\pi_n : \ell^n_{\infty} \rightarrow\ell^{I_n}_{\infty}$ is
the canonical projection.
\end{lemma}

Before coming to the statement of the main result of this section we introduce
and recall some terminology and notation. By saying that a set $\Delta$ in a
Banach space is {\it equivalent} to the standard basis $E$ of $\ell^I_1$ for
some index set $I$ we mean that there is a bijection $E\to\Delta$ which extends
to an isomorphism $\gamma : \ell^I_1 \to \overline{\spn} \Delta$. We also say
{\it $K$-equivalent} if $\| \gamma \| \| \gamma^{-1} \| \leq K$.

\begin{definition}
Let $X$ be a Banach space and $\alpha\in\IA (X)$. Let $x\in X$. We say that a
subset $I\subseteq\mathbb{Z}$ is an {\it $\ell_1$-isomorphism set for
$x$} if $\{ \alpha^i (x) : i\in I \}$ is equivalent to the standard basis of
$\ell^I_1$.
\end{definition}

\noindent Given an isometric automorphism $\alpha$ of a Banach space $X$, we
denote by $T_\alpha$ the weak$^*$ homeomorphism of the unit ball $B_1 (X^* )$
of the dual of $X$ given by $T_\alpha (\omega ) = \omega\circ\alpha$. Recall
that the {\it upper density} of a set $I\subseteq\mathbb{Z}$ is defined as
$$ \limsup_{n\to\infty} \frac{| I \cap \{ -n, -n+1 , \dots , n \} |}{2n+1} , $$
and if these ratios converge then the limit is referred to as the {\it density}
of $I$.

\begin{theorem}\label{T-zero}
Let $X$ be a Banach space and $\alpha\in\IA (X)$.
Let $Z$ be a closed $T_\alpha$-invariant subset
of $B_1(X^*)$ such that the natural linear map $X\to C(Z)$
is an isomorphism from $X$ to a (closed) linear subspace of $C(Z)$.
Then the following are
equivalent:
\begin{enumerate}
\item $\hc (\alpha ) > 0$,

\item $\htopol (T_\alpha ) > 0$,

\item $\htopol (T_\alpha ) = \infty$,

\item there exist an $x\in X$,
constants $K\geq 1$ and $d>0$, a sequence $\{n_k\}_{k\in \mathbb{N}}$ in
$\mathbb{N}$ tending to infinity, and sets $I_k\subseteq \{0, 1,
\dots, n_k-1\}$ of cardinality at least $dn_k$ such that
$\{ \alpha^i (x) : i\in I_k\}$ is $K$-equivalent to the standard basis of
$\ell^{I_k}_1$ for each $k\in\mathbb{N}$,

\item there exists an $x\in X$ with an $\ell_1$-isomorphism set of
positive density,

\item $\htopol (T_\alpha |_Z) > 0$.

\end{enumerate}
We may moreover take $x$ in (4) and (5) to be in any given total
subset $\Delta$ of $X$.
\end{theorem}

\begin{proof}
(1)$\Rightarrow$(2). The $^*$-automorphism of $C(B_1 (X^* ))$ induced from
$T_\alpha$ has positive CA entropy in view of the canonical equivariant
isometric embedding $X\hookrightarrow C(B_1 (X^* ))$ and monotonicity.
Thus $\htopol (T_\alpha ) > 0$ by Proposition~\ref{P-topol}.

(2)$\Rightarrow$(4). Without loss of generality we may assume that
$\Delta\subseteq B_1 (X)$. By assumption there exist a continuous
pseudometric $d$ on $B_1 (X^* )$, $\varepsilon > 0$, $\lambda > 0$, a sequence
$\{ n_k \}_{k\in\mathbb{N}}$ in $\mathbb{N}$ tending to infinity, and an
$(n_k , U_{d,\varepsilon} )$-separated set $E_k \subseteq B_1 (X^* )$ with
cardinality at least $e^{\lambda n_k}$ for each $k\in\mathbb{N}$. By a
simple approximation argument we may assume that there is a finite subset
$\Omega$ of $\Delta$ such that
$$ d(\sigma ,\omega ) = \sup_{x\in\Omega} | \sigma (x) - \omega (x) | $$
for all $\sigma ,\omega\in B_1 (X^* )$.

Define a linear map $\phi : X^* \to \ell^{\Omega\times n_k}_{\infty}$ by
$(\phi (f))_{x,i} = f(\alpha^i (x ))$, where the standard basis of
$\ell^{\Omega\times n_k}_{\infty}$ is indexed by
$\Omega \times \{ 0, \dots, n_k - 1 \}$. Then $\phi$ is a
contraction and $\phi (E_k )$ is $\varepsilon$-separated. By
Lemma~\ref{L-proj} there exist $d>0$ and $\delta > 0$
depending only on $\varepsilon$ and $\lambda$ such that for every
$k\geq 1$ there is a set $J_k \subseteq \Omega \times
\{ 0, \dots, n_k - 1 \}$ with
\begin{enumerate}
\item[(i)] $| J_k |\geq d|\Omega |n_k$, and

\item[(ii)] $\pi (\phi (B_1 (X^*)))\supseteq B_{\delta} (\ell^{J_k}_{\infty})$,
where $\pi : \ell^{\Omega\times n_k}_{\infty}\to \ell^{J_k}_{\infty}$ is the
canonical projection.
\end{enumerate}
Then for every $k\geq 1$ there exist some $x_k\in
\Omega$ and a set $I_k\subseteq \{ 0, \dots, n_k - 1 \}$ such that
$|I_k |\geq d n_k$ and $\{ x_k \} \times I_k \subseteq J_k$.
Consequently $\pi' (\phi (B_1 (X^* )))\supseteq B_{\delta}
(\ell^{I_k}_{\infty})$, where $\pi' : \ell^{\Omega\times n_k}_{\infty}\to
\ell^{I_k}_{\infty}$ is the canonical projection.
The dual $(\pi' \circ \phi )^*$ is an injection of
$(\ell^{I_k}_{\infty})^* = \ell^{I_k}_1$ into $X^{**}$ and the norm of the
inverse of this injection is bounded above by $\delta^{-1}$.  Notice that
$X \subseteq X^{**}$, and from our definition of $\phi$ it is clear that
$(\pi' \circ \phi )^*$ sends the standard basis element of $\ell^{I_k}_1$
associated with $i\in I_k$ to $\alpha^i (x_k )$.

Since $\Omega$ is a finite set, there is an $x\in \Omega$ such that $x_k = x$
for infinitely many $k$. By taking a subsequence of
$\{n_k \}_{k\in \mathbb{N}}$ if necessary we may assume that $x_k = x$ for all
$k\in\mathbb{N}$, and so we obtain (4).

(5)$\Rightarrow$(1). This follows from Lemma~\ref{L-ell1}.

(4)$\Rightarrow$(3). Multiplying $x$ by a scalar we may assume that
$\| x \| = 1$. On $B_1 (X^* )$ define the weak$^*$ continuous pseudometric
$$ d(\sigma , \omega ) = | \sigma (x) - \omega (x) | . $$
Denote $\spn \{\alpha^i (x) : i\in I_k \}$ by $V_k$, and
let $\gamma_k$ denote the linear map from $\ell^{I_k}_1$ to $V_k$ sending
the standard basis element of $\ell^{I_k}_1$ associated with $i\in I_k$
to $\alpha^i (x)$. For each $f\in
(\ell^{I_k}_1 )^*$ we have $(\gamma^{-1}_k )^* (f)\in V_k^*$ and
$\| (\gamma^{-1}_k )^* (f) \| \leq K \| f \|$. By the Hahn-Banach
theorem we may extend $(\gamma^{-1}_k )^* (f)$ to an element in
$X^*$ of norm at most $K \| f \|$, which we will still denote by
$(\gamma^{-1}_k )^* (f)$. Let $0 < \varepsilon < (2K)^{-1}$, and let
$M = \lfloor (2K\varepsilon)^{-1} \rfloor$ be the largest integer no greater
than $(2K\varepsilon )^{-1}$. Let $\{ g_i : i\in I_k \}$ be the standard basis
of $(\ell^{I_k}_1 )^* = \ell^{I_k}_{\infty}$. For each $f\in \{ 1,
\dots, M \}^{I_k}$ set $\tilde{f} = \sum_{i\in I_k} 2f(i)\varepsilon
g_i$. Then $f' := (\gamma^{-1}_k )^* (\tilde{f} )$ is in $B_1 (X^* )$.

We claim that the set $\{ f' : f\in \{1, \cdots, M \}^{I_k} \}$ is $(n_k,
U_{d,\varepsilon} )$-separated. Suppose $f, g \in \{1, \cdots, M\}^{I_k}$ and
$f(i) < g(i)$ for some $i\in I_k$. Then
\begin{align*}
d(T^i_{\alpha} (f'), T^i_{\alpha} (g'))
&= | (T^i_{\alpha} (f'))(x) - (T^i_{\alpha} (g'))(x) | \\
&= | f'(\alpha^i (x)) - g'(\alpha^i (x)) | \\
&= 2(g(i) - f(i))\varepsilon \\
&> \varepsilon ,
\end{align*}
establishing our claim. Therefore $\sep_{n_k} (T_{\alpha} ,
\varepsilon )\geq M^{|I_k|} \geq M^{dn_k}$. It follows that
$\htopol (T_{\alpha} )\geq d\log M$. Letting $\varepsilon\to 0$
we get $\htopol (T_{\alpha} ) = \infty$.

(3)$\Rightarrow$(2). Trivial.

(4)$\Rightarrow$(5). Let $Y_x$ be the collection of sets
$I\subseteq\mathbb{Z}$ such that the linear map from $\ell^I_1$ to
$\overline{\spn} \{ \alpha^i (x) : i\in I \}$ which sends the
standard basis element of $\ell^I_1$ associated with $i\in I$ to
$\alpha^i (x)$ is a $K$-isomorphism. If we identify subsets of
$\mathbb{Z}$ with elements of $\{ 0,1 \}^{\mathbb{Z}}$ via their
characteristic functions then $Y_x$ is a closed shift-invariant
subset of $\{ 0,1 \}^{\mathbb{Z}}$. It follows by the argument in
the second paragraph of the proof of Theorem 3.2 in \cite{GW} that
$Y_x$ has an element $J$ with density at least $d$. Clearly $J$ is
an $\ell_1$-isomorphism set for $x$.

(6)$\Rightarrow$(2). Trivial.

(5)$\Rightarrow$(6). By (5) we can find an $x\in X$ with an
$\ell_1$-isomorphism set of positive density.
Then the image of $x$ under the natural map $X\to C(Z)$
has the same $\ell_1$ isomorphism set of positive density with respect to the
induced automorphism $\alpha'$ of $C(Z)$. Since (5) implies (1) we get
$\hc (\alpha' )>0$. Then (6) follows from Proposition~\ref{P-topol}.
\end{proof}

\begin{remark}\label{R-zero}
It is easily seen that for the subset $\Delta$
we need in fact only assume that $\bigcup_{j\in\mathbb{Z}} \alpha^j (\Delta )$
is total in $X$. Furthermore, if $\alpha$ is a $^*$-automorphism of a unital
commutative $C^*$-algebra then it suffices that $\bigcup_{j\in\mathbb{Z}}
\alpha^j (\Delta )$ generate the $C^*$-algebra, since in the proof of
(2)$\Rightarrow$(4) the sets $E_k$ may be taken to lie in the pure
state space (apply (2)$\Rightarrow$(1) and Proposition~\ref{P-topol}) and
in this case pure states are multiplicative, so that we may choose the set
$\Omega$ to lie in a given $\Delta$ of the desired type.
\end{remark}

\begin{remark}\label{R-VB}
Combining Theorem~\ref{T-zero} with Lemma~\ref{L-ell1mat} yields an alternate
proof of Theorem~3.3 in \cite{EID}.
\end{remark}

\begin{corollary}\label{C-zeroisom}
Let $X_1$ and $X_2$ be Banach spaces and let $\alpha_1 \in \IA (X_1 )$ and
$\alpha_2 \in \IA (X_2 )$. Suppose there is an isomorphism
$\gamma : X_1 \to X_2$ such that
$\gamma\circ\alpha_1 = \alpha_2 \circ\gamma$. Then $\hc (\alpha_1 ) = 0$ if
and only if $\hc (\alpha_2 ) = 0$.
\end{corollary}

Since a Banach space with separable dual contains no isomorphic copy
of $\ell_1$, we have:

\begin{corollary}\label{C-sepdual}
Every isometric automorphism of a Banach space with separable dual has zero
CA entropy.
\end{corollary}

\noindent It follows, for example, that $\ell_p$ for $1<p<\infty$, $c_0$, and
the compact operators on a separable Hilbert space admit no isometric
automorphisms with nonzero CA entropy. This is also true for the
nonseparable versions of these spaces since the closed subspace dynamically
generated by a finite subset lies within a copy of the corresponding separable
version.

Since topological entropy is nonincreasing under passing to
closed invariant sets, by
Theorem~\ref{T-zero} zero CA entropy is preserved under taking quotients:

\begin{corollary}\label{C-quotient}
Let $X$ be a Banach space, $Y\subseteq X$ a closed subspace, and
$Q : X\to X/Y$ the quotient map. Let $\alpha$ be an isometric automorphism
of $X$ such that $\alpha (Y) = Y$, and suppose that $\alpha$ has zero
CA entropy. Then the induced isometric automorphism of $X/Y$ has zero
CA entropy.
\end{corollary}

The following corollary extends \cite[Thm.\ A]{GW} to the nonmetrizable case.

\begin{corollary}\label{C-nonmetric}
Let $K$ be a compact Hausdorff space and $T : K\to K$ a homeomorphism. Let 
$S_T$ be the induced homeomorphism of the space of probability
measures on $K$ with the weak$^*$ topology. Then $\htopol (T)=0$ if and only
if $\htopol (S_T )=0$.
\end{corollary}

\begin{corollary}\label{C-metrizable}
Let $K$ be a compact Hausdorff space and $T : K\to K$ a homeomorphism.
Then $\htopol (T)=0$ if and only if $\htopol (T')=0$ for
all metrizable factors $T'$ of $T$.
\end{corollary}

\begin{proof}
The ``only if'' part is immediate as topological entropy is nonincreasing under
taking factors. The ``if'' part follows from Theorem~\ref{T-zero},
Proposition~\ref{P-topol}, and the
fact that there is a one-to-one correspondence between $\alpha_T$-invariant
unital separable $C^*$-subalgebras of $C(K)$ and metrizable factors of $T$,
where $\alpha_T$ is the induced automorphism of $C(K)$.
\end{proof}

\begin{corollary}\label{C-zerosum}
Let $I$ be a nonempty index set and for each $i\in I$ let $\alpha_i$ be an
isometric automorphism of a Banach space $X_i$ with $\hc (\alpha_i ) = 0$.
Let $1\leq p < \infty$, and consider the isometric automorphism $\alpha$ of
the $\ell_p$-direct sum $\left(\bigoplus_{i\in I} X_i \right)_p$ given by
$\alpha ((x_i )_{i\in I} ) = (\alpha_i (x_i ))_{i\in I}$. Then
$\hc (\alpha ) = 0$.
\end{corollary}

\begin{proof}
By Proposition~\ref{P-prop}(ii) we may assume that $I$ is finite. Then the
formal identity map from $\left(\bigoplus_{i\in I} X_i \right)_p$ to the
$\ell_\infty$-direct sum $\left(\bigoplus_{i\in I} X_i \right)_\infty$ is a
linear contraction with bounded inverse (of norm $|I|^{1/p}$), and so we
obtain the desired conclusion from Proposition~\ref{P-directsum} and
Corollary~\ref{C-zeroisom}.
\end{proof}

\begin{corollary}\label{C-Bochner}
Let $(S, \Sigma, \mu)$ be a $\sigma$-finite measure space
and let $\alpha$ be an isometry of a Banach space $X$ with $\hc (\alpha )=0$.
Let $1\leq p<\infty$, and consider the isometric automorphism $\beta$ of
the space $L_p (S, \Sigma, \mu, X)$ given by $(\beta (f))(s)=\alpha (f(s))$.
Then $\hc (\beta )=0$.
\end{corollary}

\begin{proof}
Recall that $L_p (S, \Sigma, \mu, X)$ is the Banach space
of all (equivalence classes of) strongly $\mu$-measurable functions $f$
from $S$ into $X$ such that the norm
$$ \| f \|_p := \Big( \int_{S}\| f(s) \|^p \, d\mu(s) \Big)^{1/p} $$
is finite \cite[Chapter V]{Yos} \cite[Chapter II]{DU}.
The set of functions of the form $f=x\chi_E$, where $x\in X$ and
$\chi_E$ is the characteristic function of some $E \in \Sigma$ with
$\mu(E)<\infty$, is total in $L_p (S, \Sigma, \mu, X)$, and so the desired
conclusion follows from Corollary~\ref{C-zerosum}.
\end{proof}

In view of Theorem~\ref{T-zero} we might ask whether there exists a Banach
space which contains $\ell_1$ isometrically but does not admit an isometric
automorphism with positive CA entropy. We claim that $C(\mathbb{T} )$ is such
an example. By a theorem of Banach and Mazur $C(\mathbb{T} )$ contains
$\ell_1$ isometrically (this result is usually stated for the Cantor set
$\Delta$ or the unit interval, but we can embed $C(\Delta )$ into
$C(\mathbb{T} )$ isometrically by viewing $\Delta$ as a subset of $\mathbb{T}$
and extending functions linearly on the complement). It is also well known
that $\mathbb{T}$ admits no homeomorphisms with positive topological
entropy (see \cite{Wal}). Thus we need only appeal to the following
observation and proposition.

Let $K$ be a compact Hausdorff space and $\alpha$ an isometric automorphism of
$C(K)$. By the Banach-Stone theorem there are a unitary $w\in C(K)$ and a
homeomorphism $T: K\to K$ such that $\alpha (f) = w(f\circ T)$ for all
$f\in C(K)$.

\begin{proposition}\label{P-cms}
Let $K$, $\alpha$, $w$, and $T$ be as above. Then
$$ \hc (\alpha ) \leq \htopol (T) . $$
Also, if $\hc (\alpha ) = 0$ then $\htopol (T) = 0$.
\end{proposition}

\begin{proof}
We denote by $E_K$ the weak$^*$ compact set of extreme points of
$B_1 (C(K)^* )$, whose elements are of the form $\lambda\delta_t$ where
$\lambda\in\mathbb{T}$ and $\delta_t$ is the point mass at some $t\in K$. Let
$\pi : E_K \to K$ be the continuous surjection $\lambda\delta_t \mapsto t$.
Denoting by $T'$ the homeomorphism of $E_K$ given by
$T' (\omega ) = \omega\circ\alpha$, we then have a commutative diagram
\begin{gather*}
\xymatrix{
E_K \ar[r]^-{T'} \ar[d]_-{\pi} & E_K \ar[d]^-{\pi} \\
K \ar[r]^-{T} & K .}
\end{gather*}
Given a continuous pseudometric $d$ on $K$, we define a
continuous pseudometric $e$ on $E_K$ by $e(\eta\delta_s ,
\lambda\delta_t ) = \max (| \eta - \lambda | , d(s,t))$ for all
$\eta , \lambda\in\mathbb{T}$ and $s,t\in K$. Then for a given
$t\in K$ the map $T'$ is isometric on $\pi^{-1} (t)$ with respect
to $e$. Since the collection of pseudometrics $e$ arising as
above generate the weak$^*$ topology on $E_K$ we infer that
$\htopol (T' , \pi^{-1} (t)) = 0$. Applying Theorem~17 of \cite{Bow}
(see Lemma~\ref{L-ctssurj})
and the fact that topological entropy does not increase under
taking factors, we thus have
$$ \htopol (T) \leq \htopol (T' ) \leq \htopol (T) +
\sup_{t\in K} \htopol (T' , \pi^{-1} (t)) = \htopol (T) $$
so that $\htopol (T' ) = \htopol (T)$.

Let $\iota : C(K) \to C(E_K )$ be the isometric embedding given by
$\iota (f)(\omega ) = \omega (f)$ for all $f\in C(K)$ and $\omega\in E_K$,
and $\beta$ the $^*$-automorphism of $C(E_K )$ given by
$\beta (g)(\omega ) = g (T' \omega)$ for all $g\in C(E_K )$ and
$\omega\in E_K$. Then $\alpha$ may be viewed
as the restriction of $\beta$ to $\iota (C(K))$.
By monotonicity, Proposition~\ref{P-topol}, and the previous paragraph we
thus obtain
$$ \hc (\alpha ) \leq \hc (\beta ) = \htopol (T' ) = \htopol (T) . $$

The second assertion of the theorem follows from Theorem~\ref{T-zero} in view
of the fact that $T'$ is a restriction of the induced homeomorphism
$T_\alpha$ of $B_1 (C(K)^* )$.
\end{proof}

\section{Comparisons with matricial approximation entropies}\label{S-mat}

We briefly examine here the relation between CA entropy and its matricial
analogue for exact operator spaces, which we call CCA entropy
(completely contractive approximation entropy). The latter was introduced
for $^*$-automorphisms of exact $C^*$-algebras in \cite{PS} in which case it
was shown to coincide with Voiculescu-Brown entropy \cite[Thm.\ 3.7]{PS}.
For general references on operator spaces see \cite{ER,Pis}.

Let $X$ and $Y$ be operator spaces and $\gamma : X\to Y$ a bounded linear map.
For each $\Omega\in\Fin (X)$ and $\delta > 0$ we denote by
$\CCA (\gamma , \Omega , \delta )$ the collection of triples
$(\phi , \psi ,B)$ where
$B$ is a finite-dimensional $C^*$-algebra and $\phi : X\to B$ and
$\psi : B \to Y$ are completely contractive linear maps such that
$$ \| \psi\circ\phi (x) - \gamma (x) \| < \delta $$
for all $x\in\Omega$. By a {\it nuclear embedding} of an operator space $X$ we
mean a completely isometric linear map $\iota$ from $X$ to an operator
space $Y$ such that $\CCA (\iota , \Omega , \delta )$
is nonempty for every $\Omega\in\Fin (X)$ and $\delta > 0$. An operator space
is exact (in the sense of \cite{ER}) if and only if it admits a nuclear
embedding \cite{Kir,EOR}.

Let $X$ be an exact operator space and $\iota : X\to Y$ a
nuclear embedding. For each $\Omega\in\Fin (X)$
and $\delta > 0$ we set
$$ \rcc (\Omega , \delta ) = \inf \{ \rank\, B : (\phi , \psi , B)\in
\CA (\iota , \Omega , \delta ) \} . $$
As in the contractive approximation setting this quantity is independent
of the nuclear embedding, which can be seen in the same way using the
operator injectivity of finite-dimensional $C^*$-algebras (i.e., Wittstock's
extension theorem).

Let $\alpha$ be a completely isometric automorphism of $X$. We set
\begin{align*}
\hcc (\alpha , \Omega , \delta ) &= \limsup_{n\to\infty} \frac1n 
\log\rcc (\Omega\cup\alpha\Omega \cup\cdots\cup
\alpha^{n-1} \Omega , \delta ) , \\
\hcc (\alpha , \Omega ) &= \sup_{\delta > 0}
\hcc (\alpha , \Omega , \delta ) , \\
\hcc (\alpha ) &= \sup_{\Omega\in\Fin (X)} \hcc (\alpha , \Omega )
\end{align*}
and refer to the last quantity as the {\it completely contractive approximation
entropy} or simply {\it CCA entropy} of $\alpha$. As mentioned above, this
coincides with Voiculescu-Brown entropy when $X$ is an exact $C^*$-algebra
\cite[Thm.\ 3.7]{PS}.

The following is the matricial version of Lemma~\ref{L-ell1},
obtained by applying Lemma 3.1 of \cite{EID}.

\begin{lemma}\label{L-ell1mat}
Let $X$ be an exact operator space. Let $x_1 , \dots , x_n \in X$ and suppose 
that the linear map $\gamma : \ell_1^n \to X$ sending the $i$th standard basis
element of $\ell_1^n$ to $x_i$ for each $i=1, \dots ,n$ is an isomorphism.
Let $\delta > 0$ be such that $\delta < \| \gamma^{-1} \|^{-1}$. Then
$$ \log\rcc (\Omega , \delta ) \geq na \| \gamma \|^{-2}
(\| \gamma^{-1} \|^{-1} - \delta )^2 $$
where $a>0$ is a universal constant.
\end{lemma}

\begin{proposition}\label{P-VoiBr}
Let $\alpha$ be a completely isometric automorphism of an exact
operator space $X$. Suppose that $\hc (\alpha ) > 0$. Then
$\hcc (\alpha ) > 0$.
\end{proposition}

\begin{proof}
Apply Theorem~\ref{T-zero} and Lemma~\ref{L-ell1mat}.
\end{proof}

In general there is no inequality relating CA and CCA entropy. For instance,
the tensor product shift on $M_p^{\otimes\mathbb{Z}}$ has infinite CA entropy
(Theorem~\ref{T-shift}) but CCA entropy equal to $\log p$
\cite[Prop.\ 4.7]{Voi}, while the following example shows that it is possible
to have zero CA entropy in conjunction with positive CCA entropy.

\begin{example}\label{E-CAR}
The CAR algebra $A$ may be described as the universal unital $C^*$-algebra
generated by self-adjoint unitaries $\{ u_k \}_{k\in\mathbb{Z}}$ subject to
the anticommutation relations
$$ u_i u_j + u_j u_i = 0 $$
for all $i,j\in\mathbb{Z}$ with $i\neq j$. We define a $^*$-automorphism
$\alpha$ of $A$ by setting $\alpha (u_k ) = u_{k+1}$ for all
$k\in\mathbb{Z}$. This is an example of a bitstream $C^*$-dynamical system as
studied in \cite{GS} (see also \cite{PP}). Let
$W$ be the operator space obtained as closure of the span of
$\{ u_k \}_{k\in\mathbb{Z}}$ in $A$. Then $\alpha$ restricts to a completely
isometric automorphism of $W$, which we will denote by $\beta$.

\begin{proposition}\label{P-CAR}
We have $\hcc (\beta ) > 0$ while $\hc (\beta) =  0$.
\end{proposition}

\begin{proof}
For any $n\in\mathbb{N}$ the unitaries $u_1 , \dots , u_n$ can be realized as
tensor products of Pauli matrices (see \cite[Sect.\ 2]{Har} and
\cite[Sect.\ 9.3]{Pis}) from which it can be seen that the subset
$\{ u_1 \otimes u_1 , \dots , u_n \otimes u_n \}$ of $A\otimes A$ is
isometrically equivalent to the standard basis of $\ell^n_1$
over $\mathbb{R}$ and hence $2$-equivalent to the standard basis of $\ell^n_1$
over $\mathbb{C}$. It follows by Lemma~\ref{L-ell1mat} that
$\hcc (\alpha\otimes\alpha , \{ u_1 \otimes u_1 \} ) > 0$. Since
$\hcc (\alpha\otimes\alpha , \{ u_1 \otimes u_1 \} ) \leq
2\hspace*{0.3mm}\hcc (\alpha , \{ u_1 \} )$ by the
local tensor product subadditivity of CCA entropy (cf.\ the
proof of Proposition 3.10 in \cite{Voi}),
we obtain $\hcc (\beta ) \geq \hcc (\beta , \{ u_1 \} ) =
\hcc (\alpha , \{ u_1 \} ) > 0$.

Next we recall that the set $\{ u_k \}_{k\in I}$ is equivalent to the standard
basis of $\ell_2$ (cf.\ \cite[Sect.\ 9.3]{Pis}). Indeed given a
finite set $F\subseteq\mathbb{Z}$ and real numbers $c_k$ for $k\in F$,
the anticommutation relations between the $u_i$'s yield
$\big( \sum_{k\in F} c_k u_k \big)^2 = \sum_{k\in F} c^2_k \cdot 1$ and
hence
$$ \bigg\| \sum_{k\in F} c_k u_k \bigg\|^2 = \bigg\| \bigg(
\sum_{k\in F} c_k u_k \bigg)^2 \bigg\| = \sum_{k\in F} c^2_k , $$
from which it follows that $\{ u_k \}_{k\in F}$ is isometrically equivalent
over the real numbers to the standard basis of $(\ell^F_2 )_{\mathbb{R}}$, and
hence $2$-equivalent over the complex numbers to the standard basis of
$\ell^F_2$. We conclude by Corollary~\ref{C-sepdual} that $\hc (\beta ) = 0$.
\end{proof}
\end{example}

The above example is not a $C^*$-algebra automorphism, however, and so we
ask the following question.

\begin{question}
Is there a $^*$-automorphism of an exact $C^*$-algebra for which the
Voiculescu-Brown entropy is strictly greater than the CA entropy?
\end{question}

\section{A geometric description of the topological
Pinsker algebra}\label{S-tPa}

Let $T:K\to K$ be a homeomorphism of a compact Hausdorff space. Since zero
topological entropy is preserved under taking products and subsystems, by
a standard argument (see Corollary 2.9(1) of \cite{GM}) $T$ admits a largest
factor with zero entropy, which we will refer to as the {\it topological
Pinsker factor}, following \cite{Gl}. The corresponding $C^*$-algebra will be
called the {\it topological Pinsker algebra} and denoted $P_{K,T}$. It is an
analogue of the Pinsker $\sigma$-algebra in ergodic theory.

In \cite{ZEF} F. Blanchard and Y. Lacroix constructed the topological Pinsker
factor in the metrizable setting as the quotient
system arising from the closed $T$-invariant equivalence relation
on $K$ generated by the collection of entropy pairs. Recall that an
{\it entropy pair} is a pair $(s,t)\in K\times K$ with $s\neq t$ such
that for every two-element open cover $\mathcal{U} = \{ U,V \}$ with
$s\in\interior (K\setminus U)$ and $t\in\interior (K\setminus V)$ the local
topological entropy of $T$ with respect to $\mathcal{U}$ is nonzero
\cite{Disj}.
In \cite{EID} a description of $P_{K,T}$ for metrizable $K$ was given in terms
of local Voiculescu-Brown entropy. By applying the arguments from
\cite{EID} and Section~\ref{S-zero} we will obtain in Theorem~\ref{T-tPa}
a geometric description of $P_{K,T}$ for general $K$.

For a function $f\in C(K)$ where $K$ is a compact Hausdorff space,
we denote by $d_f$ the pseudometric on $K$ given by
$$ d_f (s,t) = | f(s) - f(t) | $$
for all $s,t\in K$. For notation relating to CCA entropy see
Section~\ref{S-mat} .

\begin{proposition}\label{P-comm}
Let $K$ be a compact Hausdorff space and $T : K\to K$ a homeomorphism.
Then for any $f\in C(K)$ the following are equivalent:
\begin{enumerate}
\item $\hcc (\alpha_T , \{ f \} ) > 0$,

\item $\hc (\alpha_T , \{ f \} ) > 0$,

\item there exists an entropy pair $(s,t)\in K\times K$ with $f(s) \neq f(t)$,

\item $h_{d_f} (T) > 0$,

\item there exist $\lambda\geq 1$, $d>0$, a sequence
$\{n_k\}_{k\in \mathbb{N}}$ in $\mathbb{N}$ tending to infinity, and sets
$I_k\subseteq \{0, 1, \dots, n_k-1\}$ of cardinality at least
$dn_k$ such that $\{ \alpha_T^i (f) : i\in I_k\}$ is $\lambda$-equivalent to
the standard basis of $\ell_1$ for each $k\in\mathbb{N}$,

\item $f$ has an $\ell_1$-isomorphism set of positive density.
\end{enumerate}
\end{proposition}

\begin{proof}
By restricting to the $\alpha_T$-invariant unital
$C^*$-subalgebra of $C(K)$ generated by $f$ we may assume that $K$ is
metrizable.

(1)$\Rightarrow$(2). Since the identity map on $C(K)$ is a nuclear embedding
(see Section~\ref{S-mat}) we need only note that a contractive
linear map from an operator space into a commutative $C^*$-algebra is
automatically completely contractive \cite[Thm.\ 3.8]{Pau}.

(2)$\Rightarrow$(3)$\Rightarrow$(4). These implications follow
from the proofs of Theorem 4.3 and Lemma 4.2, respectively, in \cite{EID}.

(4)$\Rightarrow$(5)$\Rightarrow$(6). Apply the same arguments as in the 
proofs of the respective implications
(2)$\Rightarrow$(4)$\Rightarrow$(5) in Theorem~\ref{T-zero}.

(6)$\Rightarrow$(1). By assumption there exists a set $I\subseteq\mathbb{Z}$
of density greater than some $d>0$ and an isomorphism $\gamma : \ell^I_1
\to\overline{\spn} \{ \alpha_T^i (f) : i\in I \}$ sending
the standard basis element of $\ell^I_1$ associated with
$i\in I$ to $\alpha_T^i (f)$. Let $0 < \delta < \| \gamma^{-1} \|^{-1}$.
By Lemma~\ref{L-ell1mat}, for every $n\in\mathbb{N}$ we have
\begin{gather*}
\log\rcc ( \{ \alpha_T^i (f) : i\in I \cap \{ -n , -n+1 , \dots , n \} \} ,
\delta ) \hspace*{22mm} \\
\hspace*{22mm} \geq a | I \cap \{ -n , -n+1 , \dots , n \} |
\| \gamma \|^{-2} (\| \gamma^{-1} \|^{-1} - \delta )^2
\end{gather*}
for some universal constant $a>0$. Since $I$ has density greater than $d$ and
\begin{gather*} \rcc (\{ f , \alpha_T (f) , \dots , \alpha_T^{2n} (f) \} , 
\delta ) \hspace*{45mm} \\
\hspace*{25mm} = \rcc (\{ \alpha_T^{-n} (f)  , \alpha_T^{-n+1} (f) ,\dots, 
\alpha_T^n(f) \} , \delta ) , 
\end{gather*}
we infer that
$$ \log\rcc (\{ f , \alpha_T (f) , \dots , \alpha_T^{2n} (f) \} ,
\delta )\geq d(2n+1)a \| \gamma \|^{-2}
(\| \gamma^{-1} \|^{-1} - \delta )^2 $$
for all sufficiently large $n\in\mathbb{N}$. Hence
$\hcc (\alpha_T , \{ f \}) > 0$.
\end{proof}

\begin{corollary}\label{C-Aff}
Let $A$ be a unital $C^*$-algebra and $\alpha$ an automorphism of
$A$. Let $\alpha'$ be the automorphism of $C(S(A))$ given by
$\alpha' (f)(\sigma ) = f(\sigma\circ\alpha )$ for all $f\in C(S(A))$ and
$\sigma\in S(A)$. Recall that there is an order isomorphism $x\mapsto\bar{x}$
from $A$ to the affine function space $\Aff (S(A)) \subseteq C(S(A))$ given by
$\bar{x} (\sigma ) = \sigma (x)$ for all $x\in A$ and $\sigma\in S(A)$.
Then for any $x\in A$ we have that $\hc (\alpha', \{ \bar{x} \} )>0$
implies $\hc (\alpha, \{ x \} )>0$.
\end{corollary}

\begin{proof}
The map $x\mapsto\bar{x}$ is a $2$-isomorphism of Banach spaces which
conjugates $\alpha$ to $\alpha' \big| {}_{\Aff (S(A))}$, and so we obtain the
conclusion from the implication (2)$\Rightarrow$(6) in Proposition~\ref{P-comm}
and an appeal to Lemma~\ref{L-ell1}.
\end{proof}

\begin{theorem}\label{T-tPa}
Let $T:K\to K$ be a homeomorphism of a compact Hausdorff space. Then the
topological Pinsker algebra $P_{K,T}$ is equal to the set of all $f\in C(K)$
which do not have an $\ell_1$-isomorphism set of positive density.
\end{theorem}

\begin{proof}
This follows from Proposition~\ref{P-comm} and Remark~\ref{R-zero}.
\end{proof}

\begin{corollary}\label{C-pe}
The homeomorphism $T$ has positive
topological entropy if and only if there is an $f\in C(K)$ with an
$\ell_1$-isomorphism set of positive density.
\end{corollary}

\begin{corollary}\label{C-cpe}
The homeomorphism $T$ has completely positive
entropy (i.e., every nontrivial factor has positive topological
entropy \cite{FPTE}) if and only if every nonconstant $f\in C(K)$ has an
$\ell_1$-isomorphism set of positive density.
\end{corollary}

Recently in \cite{GM} E. Glasner and M. Megrelishvili established a
Bourgain-Fremlin-\linebreak Talagrand dichotomy for metrizable topological dynamical
systems according to which the enveloping semigroup either
\begin{enumerate}
\item is a separable Rosenthal compactum (and hence has cardinality at most
$2^{\aleph_0}$), or

\item contains a homeomorphic copy of the Stone-\v{C}ech compactification
$\beta\mathbb{N}$ of $\mathbb{N}$ (and hence has
cardinality $2^{2^{\aleph_0}}$).
\end{enumerate}
A topological dynamical system (i.e., a compact space with an action
of a topological group) is said to be {\it tame} if its enveloping semigroup
is separable and Fr\'{e}chet \cite{tame}, which is equivalent to (1) in the
context of the above dichotomy.
In particular, the enveloping semigroup of a tame system has cardinality
at most $2^{\aleph_0}$.
Consider now a homeomorphism $T: K\to K$ of a
compact Hausdorff space.
If the system $(K,T)$ is tame then it is {\it $\mathbb{Z}$-regular}, i.e., 
$C(K)$ does not contain a function
$f$ such that the orbit $\{ f\circ T^n \}_{n\in\mathbb{Z}}$ admits an
infinite subset equivalent to the standard basis of $\ell_1$ (the
$\mathbb{N}$-system version of this property is called {\it regularity}
in \cite{K}), for otherwise
the enveloping semigroup would contain a homeomorphic copy of
$\beta\mathbb{N}$ (see the proof of Corollary 5.4 in \cite{K}).
Thus Corollary~\ref{C-pe} yields the following, which generalizes a result
of Glasner \cite{tame}, who proved it for $(K,T)$ metrizable and minimal.

\begin{corollary}
A tame homeomorphism $T:K\to K$ of a compact Hausdorff space has zero
topological entropy.
\end{corollary}

If the system $(K,T)$ is HNS (hereditarily not sensitive) \cite[Defn.\ 9.1]{GM}
and $K$ is metrizable then it is tame and hence has zero topological entropy
by the above corollary. More generally we have:

\begin{corollary}
An HNS homeomorphism $T:K\to K$ of a compact Hausdorff space has zero
topological entropy.
\end{corollary}


\begin{proof}
By Theorems 9.8 and 7.6 of \cite{GM} the system $(K,T)$ is HNS if and only if
every $f\in C(K)$ is in $\Asp (K)$, i.e., if and only if the pseudometric
space $(K, \rho_{\mathbb{Z}, f}|_K)$ is separable, where $\rho_{\mathbb{Z}, f}$
is the pseudometric on $C(K)^*$ defined by
$$ \rho_{\mathbb{Z}, f}(\sigma, \omega)=\sup_{n\in\mathbb{Z}}
|\sigma (f\circ T^n )-\omega (f\circ T^n )| . $$
By Lemma 1.5.3 of \cite{Fa} the pseudometric space
$(C(K)^*, \rho_{\mathbb{Z}, f})$
is also separable. Thus the orbit $\{f\circ T^n\}_{n\in\mathbb{Z}}$ admits no
infinite subset equivalent to the standard basis of $\ell_1$, and we obtain
the result by Corollary~\ref{C-pe}.
\end{proof}

Versions of
Proposition~\ref{P-comm} and Theorem~\ref{T-tPa} can be similarly established
for sequence topological entropy \cite{SE,NSSEP}. Recall that a system $(K,T)$ 
consisting of a homeomorphism $T:K\to K$ of a compact Hausdorff space is said 
to be {\it null} if its sequence topological
entropy is zero for all sequences. Nullness is preseved under taking products
and subsystems, and so every system admits a largest null factor. 
In analogy with Theorem~\ref{T-tPa} we then have the following.

\begin{theorem}\label{T-null}
Let $T:K\to K$ be a homeomorphism of a compact Hausdorff space. Then the 
largest null factor of the system $(K,T)$, when viewed as a
dynamically invariant $C^*$-subalgebra of $C(K)$, is equal to the set
of all $f\in C(K)$ satisfying the property that for every $\lambda\geq 1$ there
exists an $m\in\mathbb{N}$ such that if $\Omega$ is a subset of
$\{ f\circ T^n \}_{n\in\mathbb{N}}$ which is $\lambda$-equivalent to the
standard basis of $\ell^{\Omega}_1$ then $|\Omega | \leq m$. In particular, the
system $(K,T)$ is null if and only if the above property is satisfied by
every $f\in C(K)$ (equivalently, by every $f$ in a given $\Delta\subseteq C(K)$
such that $\bigcup_{j\in\mathbb{Z}} \alpha^j (\Delta )$ generates $C(K)$ 
as a $C^*$-algebra). 
\end{theorem}

\begin{corollary}
A null homeomorphism $T:K\to K$ of a compact Hausdorff space is 
$\mathbb{Z}$-regular.
\end{corollary}

\noindent We also have the following analogue for nullness of 
\cite[Thm.\ A]{GW} (cf.\ Corollary~\ref{C-nonmetric}).

\begin{theorem}\label{T-nullstate}
A homeomorphism $T:K\to K$ of a compact Hausdorff space is null if and only if
the induced weak$^*$ homeomorphism of the space $M(K)$ of probability measures 
on $K$ is null.
\end{theorem}

\begin{proof}
The image of $C(K)$ in $C(M(K))$ under the equivariant map given by evaluation 
generates $C(M(K))$ as a $C^*$-algebra, and so we can apply 
Theorem~\ref{T-null} to obtain the nontrivial direction.
\end{proof}

\noindent We point out that, for minimal distal metrizable systems, 
nullness, $\mathbb{Z}$-regularity, and equicontinuity are equivalent. 
The equivalence of nullness and equicontinuity is established in 
\cite{NSSEP} (see Corollaries~2.1(2) and 4.2 therein), while the 
equivalence of $\mathbb{Z}$-regularity and equicontinuity follows from 
Corollary~1.8 of \cite{tame} and the fact that tameness and 
$\mathbb{Z}$-regularity coincide in the metrizable case. The 
equivalence of nullness and $\mathbb{Z}$-regularity can also be extracted from
\cite{NSSEP}: if the system is not null then the proofs of 
Corollaries~4.2 and 4.1 and Lemma~3.1 in \cite{NSSEP} show that there is 
a two-element open cover satisfying the
property in the statement of Proposition~2.3 in \cite{NSSEP}, yielding the 
existence of a real-valued continuous function on the space whose forward 
orbit contains a subsequence isometrically equivalent to the standard basis of
$\ell_1$ over $\mathbb{R}$, so that the system is not $\mathbb{Z}$-regular.

\section{A dynamical characterization of type I $C^*$-algebras}\label{S-typeI}

By \cite{NS,typeI} a separable unital $C^*$-algebra is type I if and only if
every inner $^*$-automor-\linebreak phism has zero Connes-Narnhofer-Thirring entropy 
with respect to each invariant state. The question thus arises of whether there is
a topological version of this result, and indeed N. P. Brown conjectured
in \cite{typeI} that the analogous assertion for Voiculescu-Brown entropy
also holds. One major difficulty is that the behaviour of zero
Voiculescu-Brown entropy with respect to taking extensions is not well
understood. For CA entropy, however, we can show that
zero values persist under taking extensions by reducing the problem to a
topological-dynamical one by means of Theorem~\ref{T-zero}. The
topological-dynamical ingredient that we require is
provided by the following lemma, which, using the notions of separated and
spanning sets described in Section~\ref{S-zero}, can be established in the
same way as its specialization to the metric setting \cite[Thm.\ 17]{Bow}.

\begin{lemma}\label{L-ctssurj}
Let $K,J$ be compact Hausdorff spaces and $T:K\to K$, $S:J\to J$
homeomorphisms. Let $\pi : K\to J$ be a continuous surjective map
such that $\pi\circ T = S\circ\pi$. Then
$$ \htopol (T) \leq \htopol (S) + \sup_{s\in J} \htopol
(T, \pi^{-1} (s)) . $$
\end{lemma}

\begin{lemma}\label{L-quotient}
Let $X$ be a Banach space, $Y\subseteq X$ a closed subspace, and
$Q : X\to X/Y$ the quotient map. Let $\alpha$ be an isometric automorphism
of $X$ such that $\alpha (Y) = Y$, and denote by $\bar{\alpha}$ the
induced isometric automorphism of $X/Y$. Then $\hc (\alpha)=0$ if and only if
$\hc (\alpha |_Y ) = \hc (\bar{\alpha} ) = 0$.
\end{lemma}

\begin{proof}
The ``only if'' part follows from monotonicity and Corollary~\ref{C-quotient}.
For the ``if'' part we first observe that we have a commutative diagram
\begin{gather*}
\xymatrix{
B_1 (X^* ) \ar[r]^-{T_\alpha}
\ar[d]_-{\Phi} & B_1 (X^* ) \ar[d]^-{\Phi} \\
B_1 (Y^* ) \ar[r]^-{T_{\alpha | Y}} & B_1 (Y^* )}
\end{gather*}
where $\Phi$ is the weak$^*$ continuous surjective map given by restriction.
Lemma~\ref{L-ctssurj} then yields
$$ \htopol (T_\alpha ) \leq \htopol
(T_{\alpha | Y}) + \sup_{\sigma\in B_1 (Y^* )} \htopol
(T_\alpha , \Phi^{-1}(\sigma )) . $$
By Theorem~\ref{T-zero} an isometric automorphism of a Banach space has
zero CA entropy if and only if the induced homeomorphism of the unit ball of
the dual has zero topological entropy, and thus, since $\hc (\alpha |_Y ) = 0$
by hypothesis, the proof will be complete once we show that $\htopol
(T_\alpha , \Phi^{-1}(\sigma )) = 0$ for all $\sigma\in B_1 (Y^* )$.

So let $\sigma\in B_1 (Y^* )$. Pick an $\omega\in\Phi^{-1} (\sigma )$
and define the weak$^*$ continuous map $\Psi : B_1 ((X/Y)^* ) \to X^*$
by $\Psi (\rho ) = 2\rho\circ Q + \omega$ for all $\rho\in B_1 ((X/Y)^* )$.
Let $\{ x_1 , x_2 , \dots , x_r \}$ be a finite subset of $B_1 (X)$ and define
on $X^*$ the weak$^*$ continuous pseudometric
$$ d(\eta , \tau ) = \sup_{1\leq j\leq r} | \eta (x_j ) - \tau (x_j ) | . $$
Note that $d(c\eta , c\tau ) = c d(\eta , \tau )$ and
$d(\eta + \rho , \tau + \rho ) = d(\eta , \tau )$ for all
$\eta , \tau , \rho\in X^*$ and $c>0$. We also define on $(X/Y)^*$ the
weak$^*$ continuous pseudometric
$$ \bar{d} (\eta , \tau ) = \sup_{1\leq j\leq r} | \eta (Q(x_j )) -
\tau (Q(x_j )) | . $$
Let $n\in\mathbb{N}$ and $\varepsilon > 0$, and let $E\subseteq B_1 ((X/Y)^* )$
be an $(n, U_{\bar{d},\varepsilon /4} )$-spanning set with
respect to $T_{\bar{\alpha}}$. Now suppose $\eta\in
\Phi^{-1} (\sigma )$. Since $Q$ is a quotient map there is a
$\rho\in (X/Y)^*$ such that $\rho\circ Q = (\eta - \omega )/2$ and
$\| \rho \| = \| \eta - \omega \| /2 \leq 1$. We can then find a
$\tau\in E$ such that
$$ \bar{d} (T^k_{\bar{\alpha}} (\rho ) , T^k_{\bar{\alpha}} (\tau )) <
\varepsilon /4 $$
for each $k = 0, \dots , n-1$. Then for each $k = 0, \dots , n-1$ we have
\begin{align*}
d(T^k_\alpha (\eta ) , T^k_\alpha (\Psi (\tau ))) &=
d(T^k_\alpha (2 \rho\circ Q + \omega ) , T^k_\alpha
(2 \tau\circ Q + \omega )) \\
&= 2d(T^k_\alpha (\rho\circ Q) , T^k_\alpha (\tau\circ Q)) \\
&= 2\bar{d} (T^k_{\bar{\alpha}} (\rho) , T^k_{\bar{\alpha}} (\tau)) \\
&< 2 (\varepsilon /4) = \varepsilon /2 .
\end{align*}
For every $\tau\in E$ we pick, if possible, an
$\eta_\tau \in\Phi^{-1} (\sigma )$ such that
$$ d( T^k_\alpha (\eta_\tau ) , T^k_\alpha (\Psi (\tau )))
< \varepsilon /2 $$
for each $k = 0, \dots , n-1$. For those $\tau\in E$ for which this
is not possible we set $\eta_\tau = 0$. It is then easily checked that
the set $F = \{ \eta_\tau : \tau\in E \}$, which has cardinality at most
that of $E$, is $(n, U_{d,\varepsilon} )$-spanning
with respect to $T_\alpha$. Hence $\spa_n (T_\alpha , \Phi^{-1} (\sigma ),
U_{d,\varepsilon} ) \leq \spa_n (T_{\bar{\alpha}} , 
U_{\bar{d},\varepsilon /4} )$.
Since $\htopol (T_{\bar{\alpha}}) = \hc (\bar{\alpha} ) = 0$
by Theorem~\ref{T-zero} and our hypothesis, we conclude
that $\htopol (T_\alpha , \Phi^{-1}(\sigma )) = 0$, as desired.
\end{proof}

We also need to know that zero CA entropy is well behaved with respect to
continuous fields over locally compact Hausdorff spaces.




\begin{lemma}\label{L-field}
Let $Z$ be a locally compact Hausdorff space, and let $(X_z )_{z\in
Z}$ be a continuous field of Banach spaces over $Z$.
Let $X$ be the Banach space of continuous sections of $(X_z )_{z\in Z}$
vanishing at infinity. Let $\alpha$ be an isometric automorphism of $X$ which
arises from $\alpha_z \in\IA (X_z )$ for $z\in Z$. Then
$\hc (\alpha )=0$ if and only if $\hc (\alpha_z )=0$ for all $z\in Z$.
\end{lemma}

\begin{proof} The ``only if'' part follows from Corollary~\ref{C-quotient}.
Suppose then that $\hc (\alpha_z )=0$ for all $z\in Z$, and let us show that
$\hc (\alpha )=0$. We first reduce the problem to the case $Z$ is compact.
For each compact subset $W\subseteq Z$ let $X_W$ be the Banach space of
continuous sections of the restriction field $(X_w )_{w\in W}$ of Banach
spaces over $W$ and let $Q_W : X\rightarrow X_W$ be the quotient map. The
isometric automorphisms $(\alpha_w )_{w\in W}$ give us an isometric
automorphism $\alpha_W$ of $X_W$. If $\hc (\alpha )>0$ then by
Theorem~\ref{T-zero} we can find an
$x\in X$ with an $\ell_1$-isomorphism set of positive density.
Using the fact that $x$ vanishes at infinity, we can find a
compact subset $W\subseteq Z$ such that $Q_W (x)$ has the same
$\ell_1$-isomorphism set of positive density under $\alpha_W$.
Then $\hc (\alpha_W )>0$ by Theorem~\ref{T-zero}. Replacing $Z$ by
$W$, we may thus assume that $Z$ is compact.

Next we reduce the problem to the case in which $X$ contains an
element $x$ such that $Q_z (x)$ has norm $1$ and is fixed under
$\alpha_z$ for every $z\in Z$, where $Q_z$ the quotient map $X\to
X_z$. Set $Y_z$ to be the $\ell_{\infty}$-direct sum $X_z\oplus
\mathbb{C}$ for each $z\in Z$. Then $(Y_z )_{z\in Z}$ is in a natural way
a continuous field of Banach spaces over $Z$. The global section
space $Y$ of this field is the $\ell_{\infty}$-direct sum $X\oplus
C(Z)$. The automorphisms $\alpha$ and $\alpha_z$ naturally extend
to isometric automorphisms of $Y$ and $Y_z$ fixing
$C(Z)$ and $\mathbb{C}$, respectively. By Lemma~\ref{L-quotient} we have
$\hc(\gamma_z )=0$ for every $z\in Z$, where $\gamma_z$ is the extension
of $\alpha_z$.
Note the constant function $1\in C(T)\subseteq Y$ satisfies the
above requirement. Replacing $X$ by $Y$, we may assume that $X$
contains such an element $x$.

For each $z\in Z$ denote by $S_z$ the subset of $B_1 (X^*_z )$ consisting 
of linear functionals $\sigma$ satisfying $\sigma(Q_z(x))\geq 1/2$. 
Clearly the sets $Q^*_z (S_z )$ for $z\in Z$ are pairwise disjoint. Let
$S=\bigcup_{z\in Z} Q^*_z (S_z )$, and let $\psi : S\rightarrow Z$ be the map
which sends $\sigma_z \in Q^*_z (S_z )$ to $z\in Z$.
Denote by $T_{\alpha}$ (resp.\ $T_{\alpha_z}$)
the homeomorphism of $B_1 (X^* )$ (resp.\ $B_1 (X^*_z )$) induced by $\alpha$
(resp.\ $\alpha_z$). One checks easily that $\psi$
is surjective and continuous, and that $S$ is a $T_{\alpha}$-invariant
closed subset of $B_1 (X^* )$. By
Theorem~\ref{T-zero} we have $\htopol (T_{\alpha_z} |_{S_z})=0$.
Applying Lemma~\ref{L-ctssurj} we get
$$ \htopol (T_{\alpha} |_S ) = \sup_{z\in Z} \htopol
(T_{\alpha},  \psi^{-1} (z)) = \sup_{z\in Z} \htopol
(T_{\alpha_z} |_{S_z} ) = 0 . $$
Now given any unit vectors $v$ and $w$ in a Banach space $V$ and $r\in [0,1]$,
there exists a $\sigma\in B_1 (V^* )$ with $\sigma (v) \geq r$ and 
$| \sigma (w) | \geq (1-r)/3$ (indeed if $\sigma$ is an element of 
$B_1 (V^* )$ such that $\sigma (v) = 1$ and $| \sigma (w) | < (1-r)/3$, then 
choose a $\tau\in B_1 (V^* )$ with $\tau (v) \geq 0$ and $| \tau (w) | = 1$ and
replace $\sigma$ with $(\sigma + (1-r) \tau ) / \| \sigma + (1-r) \tau \|$).
It follows that the natural linear map $X\to C(S)$ is an
isomorphism from $X$ to a closed linear subspace of $C(S)$. 
We thus conclude by Theorem~\ref{T-zero} that $\hc (\alpha )=0$.
\end{proof}

Given a $C^*$-algebra $A$ we denote by $M(A)$ its multiplier algebra and by
$\Prim (A)$ its primitive ideal space.

\begin{theorem}\label{T-typeI}
Let $A$ be a $C^*$-algebra. Then the following are equivalent:

\begin{enumerate}
\item $A$ is type I,

\item $\hc(\alpha) = 0$ for every $\alpha\in \Aut(A)$ with trivial induced
action on $\Prim (A)$,

\item $\hc(\Ad u) = 0$ for every unitary $u\in M(A)$,

\item $\hc(\Ad u) < \infty$ for every unitary $u\in M(A)$.
\end{enumerate}
\end{theorem}

\begin{proof}
(1)$\Rightarrow$(2). By \cite[Thm.\ 6.2.11]{Ped} there is a composition series
$(I_{\rho})_{0\leq\rho\leq\mu}$ (i.e., an increasing family of closed
two-sided ideals of $A$ indexed by ordinals with $I_\mu = A$ and
$I_{\rho}$ equal the norm closure of $\bigcup_{\rho' < \rho} I_{\rho'}$ for
any limit ordinal $\rho\leq\mu$) such that each quotient
$I_{\rho+1}/I_{\rho}$ for $0\leq\rho < \mu$ is a continuous trace
$C^*$-algebra. By Lemma~\ref{L-quotient} and Proposition~\ref{P-prop} we may
therefore assume that $A$ is a continuous trace $C^*$-algebra. Thus
$\hat{A}$ is a locally compact
Hausdorff space and $A$ is the $C^*$-algebra of continuous sections vanishing
at infinity of a continuous field of $C^*$-algebras over $\hat{A}$ with each
fibre $A_x$ equal to the compact operators $\mathcal{K}(\mathcal{H}_x)$ for
some Hilbert space $\mathcal{H}_x$. By Corollary~\ref{C-sepdual} (see
the comment following it) every isometric automorphism of
$\mathcal{K}(\mathcal{H}_x)$ has zero CA entropy, and so by 
Lemma~\ref{L-field} we conclude that $\hc(\alpha)=0$.

(2)$\Rightarrow$(3)$\Rightarrow$(4). Trivial.

(4)$\Rightarrow$(1). Here we can apply the construction in the proof of
Theorem~1.2 in \cite{typeI}. Suppose that $A$ is not type I. Denote by
$\tilde{A}$ the unitization of $A$. Let $\gamma$ be the tensor product shift
on the CAR algebra $M_{2^\infty} = M_2^{\otimes\mathbb{Z}}$. Set
$\beta = \gamma^{\otimes\mathbb{N}} \in
\Aut ((M_{2^\infty} )^{\otimes\mathbb{N}})$.
As in the proof of Theorem~1.2 in \cite{typeI}, we can find a unital
$C^*$-subalgebra $B\subseteq\tilde{A}$ such that
there are a unitary $u\in\tilde{A}$ and a surjective $^*$-homomorphism
$\pi : B\rightarrow (M_{2^\infty} )^{\otimes\mathbb{N}}$ with $\pi\circ\Ad u =
\beta\circ\pi$. For each $j\in\mathbb{N}$ let $x_j$ be the matrix
$\big[ \begin{smallmatrix} 1 & 0 \\ 0 & -1 \end{smallmatrix} \big]$
viewed as an element of the zeroeth copy of $M_2$ in the $j$th copy of
$M_{2^\infty}$ in $(M_{2^\infty} )^{\otimes\mathbb{N}}$, and pick a
$b_j \in B$ such that $\| b_j \| = 1$ and $\pi (b_j ) = x_j$. Let
$m\in\mathbb{N}$ and set $\Omega_m = \{ b_j : j=1, \dots , m \}$. For each
$b\in\Omega_m$ write $b = \lambda_b 1
+ a_b$ where $\lambda_b \in\mathbb{C}$ and $a_b\in A$, and note that
$\| a_b \| \leq 2$ (otherwise $a_b$ would be invertible).

Denote by $\varphi$ the contractive linear map from $\ell^{\{ 1, \dots , m \}
\times\mathbb{Z}}_1$ to $\overline{\spn} \bigcup_{k\in\mathbb{Z}}
\Ad u^k (\Omega_m )$ which sends standard basis elements to elements in
$\bigcup_{k\in\mathbb{Z}} \Ad u^k (\Omega_m )$ respecting the
indexing in the obvious way. It is easily checked that $\varphi$ is a
$2$-isomorphism. Set $\Omega'_m = \{ a_b : b\in\Omega_m \}$ and 
let $\varphi'$ be the bounded linear map from
$\ell^{\{ 1, \dots , m \} \times\mathbb{Z}}_1$ to
$\overline{\spn} \bigcup_{k\in\mathbb{Z}} \Ad u^k (\Omega_m' )$ which
corresponds to $\varphi$ on the standard basis elements via the association
of $\Ad u^k (a_b )$ with $\Ad u^k (b)$.
We will argue that $\varphi'$ has an inverse of norm at most $2$.
Suppose $f$ is a norm one element of
$\ell^{\{ 1, \dots , m \} \times \{ 1, \dots , n \}}_1$. Then
$g = f\oplus (-f) \in\ell^{\{ 1, \dots , m \} \times \{ 1, \dots , 2n \}}_1$
has norm $2$, and so $\| \varphi (g) \| \geq 1$. Since $\Ad u$ is unital we
have
$$ \varphi (g) = \varphi (f) - u^n \varphi (f) u^{-n}
= \varphi' (f) - u^n \varphi' (f) u^{-n} = \varphi' (g) $$
and hence at least one of $\| \varphi' (f) \|$ and
$\| u^n \varphi' (f) u^{-n} \|$ is greater than or equal to $1/2$.
Therefore $\| \varphi' (f) \| \geq 1/2$, as desired.

Having shown that $\varphi'$ is an isomorphism (in fact a $4$-isomorphism
since $\| \varphi' \| \leq 2$), it follows from Lemma~\ref{L-ell1} that for
a given $0 < \delta < 1/2$ we have $\hc(\Ad u, \Omega'_m, \delta)\geq am$ for
some $a>0$ which does not depend on $m$, and consequently $\hc(\Ad u)=\infty$.
\end{proof}

\section{Tensor products, crossed products, and free products}\label{S-prod}

We have tensor product subadditivity with respect to the injective tensor
product:

\begin{proposition}\label{P-tensorprod}
Let $X_1$ and $X_2$ be Banach spaces and let $\alpha_1 \in \IA (X_1 )$ and
$\alpha_2 \in \IA (X_2 )$. Then the isometric automorphism
$\alpha_1 \check{\otimes}\alpha_2$ of $X_1 \check{\otimes} X_2$ satisfies
$$ \hc (\alpha_1 \check{\otimes}\alpha_2 ) \leq
\hc (\alpha_1 ) + \hc (\alpha_2 ) . $$
\end{proposition}

\begin{proof}
Let $\iota_1 : X_1 \to Y_1$ and $\iota_2 : X_2 \to Y_2$ be CA embeddings.
Then $\iota_1 \check{\otimes} \iota_2 : X_1 \check{\otimes} X_2 \to
Y_1 \check{\otimes} Y_2$ is a CA embedding, and if for given
finite subsets $\Omega_1 \subseteq B_1 (X_1 )$ and $\Omega_2 \subseteq
B_1 (X_2 )$ and $\delta > 0$
we have $(\phi_i , \psi_i , d_i )\in\CA (\iota_i , \Omega_i , \delta )$
for $i=1,2$ then since $\ell_\infty^{d_1} \check{\otimes}
\ell_\infty^{d_2} = \ell_\infty^{d_1 d_2}$ it is readily checked that
$$ (\phi_1 \check{\otimes} \phi_2 , \psi_1 \check{\otimes} \psi_2 ,
d_1 d_2 )\in\CA (\iota_1 \check{\otimes} \iota_2 , \Omega_1 \otimes
\Omega_2 ,  2\delta ) , $$
from which the result follows.
\end{proof}

\noindent The subadditivity of Proposition~\ref{P-tensorprod} can fail for
other tensor product norms, in particular for the minimal operator space
tensor product, as the following example demonstrates. Let $\beta$ be the
completely isometric automorphism of the operator subspace $W$ of the
CAR algebra as described in Example~\ref{E-CAR}. The CAR algebra is
$^*$-isomorphic to the infinite tensor product UHF algebra
$M_2^{\otimes\mathbb{Z}}$, and so $W$ is exact.
In Proposition~\ref{P-CAR} it was shown that $\hc (\beta ) = 0$.
On the other hand, applying monotonicity and using Corollary~\ref{C-cpe} as in
the first part of the proof of Proposition~\ref{P-CAR} (to which we refer for
notation) we have
$$ \hc (\beta\otimes\beta ) \geq
\hc (\alpha\otimes\alpha , \{ u_0 \otimes u_0 \} ) > 0 . $$

We turn next to $C^*$-crossed products. All crossed products considered here
will be reduced and so for economy we won't bother to tag the product
symbol to indicate this, contrary to usual practice. Note however that
the full and reduced crossed products coincide if the acting group is
amenable, as will ultimately be the case here.

We will establish
in Theorem~\ref{T-crossedprod} below an analogue of Theorem~5.3 of
\cite{PS}, whose proof based on Imai-Takai duality we will adapt. In our
case we cannot establish equality beyond the zero entropy case
unless the action is on a commutative $C^*$-algebra, since tensor product
subadditivity only holds with respect to the injective tensor product.

The following lemma and proposition are the analogues of Lemma~2.2
and Theorem~5.2, respectively, in \cite{PS}.

\begin{lemma}\label{L-factor}
Let $X_1$ and $X_2$ be Banach spaces and $\alpha_1 \in \IA (X_1 )$ and
$\alpha_2 \in \IA (X_2 )$. Suppose there exists a net
\begin{gather*}
\xymatrix{
X_1 \ar[r]^-{S_\lambda} & X_2 \ar[r]^-{T_\lambda} & X_1}
\end{gather*}
of contractive linear maps such that $T_\lambda \circ S_\lambda$ converges
to $\id_{X_1}$ in the point-norm topology and $S_\lambda \circ\alpha_1 =
\alpha_2 \circ S_\lambda$ and $T_\lambda \circ\alpha_2 = \alpha_1 \circ
T_\lambda$ for all $\lambda$. Then
$\hc (\alpha_1 ) \leq \hc (\alpha_2 )$.
\end{lemma}

\begin{proof}
Let $\iota_1 : X_1 \to Y_1$ and $\iota_2 : X_2 \to Y_2$ be CA embeddings.
We may assume that $Y_1$ is injective by taking $\iota_1$ to be, for example,
the composition of the map $X_1 \to C(B_1 (X^* ))$ defined via evaluation with
the canonical embedding of $C(B_1 (X^* ))$ into its second dual.
Let $\Omega\in\Fin (X_1 )$ and $\delta > 0$. Pick a $\lambda$ such that
$\| T_\lambda \circ S_\lambda (x) - x \| < \delta$ for all $x\in\Omega$.
Let $n\in\mathbb{N}$ and suppose $(\phi , \psi , d)\in\CA
(\iota_2 , S_\lambda (\Omega ) \cup \alpha_2 (S_\lambda (\Omega )) \cup\cdots
\cup\alpha_2^{n-1} (S_\lambda (\Omega )) , \delta )$.
By the injectivity of $Y_1$ we can extend
$T_\lambda \circ\iota_2^{-1} |_{\iota_2 (X_2 )}$ to
a contractive linear map $\gamma : Y_2 \to Y_1$. From our
assumption we have $S_\lambda \circ\alpha^k_1 (x) =
\alpha^k_2 (S_\lambda (x))$ and $T_\lambda \circ\alpha^k_2 (S_\lambda (x))
= \alpha^k_1 \circ T_\lambda (S_\lambda (x))$ for all $x\in\Omega$ and
$k\in\mathbb{Z}$, and so by an estimate using the triangle inequality we have
$$ (\phi\circ S_\lambda , \tau\circ\gamma\circ\psi , d) \in
\CA (\iota_1 , \Omega\cup \alpha_1 \Omega \cup\cdots
\cup\alpha_1^{n-1} \Omega , 2\delta ) . $$
We infer that $\hc (\alpha_1 , \Omega , 2\delta ) \leq
\hc (\alpha_2 , S_\lambda (\Omega ) , \delta )$, from which we conclude
that $\hc (\alpha_1 ) \leq \hc (\alpha_2 )$.
\end{proof}

\begin{proposition}\label{P-crossedprodleq}
Let $A$ be a $C^*$-algebra, $G$ a locally compact group, and $\alpha$ a
strongly continuous action of $G$ on $A$ by $^*$-automorphisms.
Let $g\in G$, and suppose that $e$ has a basis of neighbourhoods $N$
satisfying $gNg^{-1} = N$. Then
$$ \hc (\alpha_g ) \leq \hc (\Ad \lambda_g |_{A\rtimes_{\alpha} G} ) . $$
\end{proposition}

\begin{proof}
In the proof of Theorem~5.2 in \cite{PS} it is shown that the pair
$$ (A, \alpha_g ), (A\rtimes_\alpha G,
\Ad \lambda_g |_{A\rtimes_{\alpha} G} ) $$
has the covariant completely contractive factorization property
(see Section~2 of \cite{PS}), and so we can apply
Lemma~\ref{L-factor} to obtain the result.
\end{proof}



For a locally compact group $G$ we consider $\Aut (G)$
with its standard topological group structure, for which a neighbourhood
basis of $\id_G$ is formed by the sets
$\{ \alpha\in\Aut (G) : \alpha (x) \in Vx \text{ and } \alpha^{-1} (x)
\in Vx \text{ for all } x\in K \}$ where $K\subseteq G$ is compact and $V$ is a
neighbourhood of $e$ in $G$.

\begin{theorem}\label{T-crossedprod}
Let $A$ be a $C^*$-algebra, $G$ an amenable locally compact
group, and $\alpha$ a strongly continuous
action of $G$ on $A$ by $^*$-automorphisms. Let $g\in G$, and suppose that
the closure of the subgroup of $\Aut (G)$ generated by $\Ad g$ is compact.
Then $\hc (\Ad \lambda_g |_{A\rtimes_{\alpha} G} ) = 0$ if and only if
$\hc (\alpha_g ) = 0$. If $A$ is furthermore assumed to be commutative then
$\hc (\Ad \lambda_g |_{A\rtimes_{\alpha} G} ) = \hc (\alpha_g )$.
\end{theorem}

\begin{proof}
By \cite[Lemma 5.1]{PS} and Proposition~\ref{P-crossedprodleq} it suffices to
prove that $\hc (\alpha_g ) = 0$ implies
$\hc (\Ad \lambda_g |_{A\rtimes_{\alpha} G} )
= 0$, and, in the case that $A$ is commutative, that
$\hc (\Ad \lambda_g |_{A\rtimes_{\alpha} G} ) \leq \hc (\alpha_g )$.

It is shown in the proof of Theorem~5.3 in \cite{PS} that the pair
$$ (A\rtimes_{\alpha} G ,
\Ad \lambda_g |_{A\rtimes_{\alpha} G} ), (A\rtimes_{\alpha} G
\rtimes_{\hat{\alpha}} G , \Ad (\lambda_g \otimes l_g r_g )
|_{A\rtimes_{\alpha} G \rtimes_{\hat{\alpha}} G} ) $$
has the covariant completely contractive factorization property
(see Section~2 of \cite{PS}) and hence by Lemma~\ref{L-factor} we have
$$ \hc (\alpha_g ) \leq \hc (\Ad \lambda_g |_{A\rtimes_{\alpha} G} ) \leq
\hc (\Ad (\lambda_g \otimes l_g r_g )
|_{(A\rtimes_{\alpha} G) \rtimes_{\hat{\alpha}} G} ) . $$
By Theorem~4.1 of \cite{PS} we have
$$ \Ad (\lambda_g \otimes l_g r_g )
|_{A\rtimes_{\alpha} G \rtimes_{\hat{\alpha}} G} =
\alpha_g \otimes (\Ad l_g r_g |_{\mathcal{K} (L^2 (G))} ) . $$
Since the minimal $C^*$-tensor product coincides with the injective tensor
product if one of the factors is commutative, it follows by
Proposition~\ref{P-tensorprod} that if $A$ is commutative then
$$ \hc(\Ad (\lambda_g \otimes l_g r_g )
|_{A\rtimes_{\alpha} G \rtimes_{\hat{\alpha}} G} )
\leq \hc (\alpha_g ) + \hc (\Ad l_g r_g |_{\mathcal{K} (L^2 (G))} ) $$
and hence $\hc (\Ad \lambda_g |_{A\rtimes_{\alpha} G} )
= \hc (\alpha_g )$ since every $^*$-automorphism of
the compact operators has zero CA entropy by Corollary~\ref{C-sepdual}.

Assuming now that $A$ is not commutative, we suppose that $\hc
(\alpha_g ) = 0$. For economy we set $\beta_g = \Ad l_g r_g
|_{\mathcal{K} (L^2 (G))}$. Let $x\in A\otimes\mathcal{K} (L^2
(G))$. To obtain the desired equality $\hc (\Ad \lambda_g
|_{A\rtimes_{\alpha} G} ) = 0$, it suffices by the general
observations above and Theorem~\ref{T-zero} to show that $x$
admits no $\ell_1$-isomorphism set of positive density with
respect to the $^*$-automorphism $\alpha_g \otimes \beta_g$. Since
the span of rank one projections is dense in $\mathcal{K} (L^2
(G))$, we may assume that $x$ is of the form $\sum_{k=1}^r a_k
\otimes p_k$ where $p_1 , \dots , p_r$ are rank one projections in
$\mathcal{B} (L^2 (G))$. Denoting by $p_\xi$ the orthogonal
projection onto the subspace spanned by a given vector $\xi\in L^2
(G)$, it is readily verified that $\beta_g (p_\xi ) =
p_{\xi\circ\Ad g^{-1}}$ for all $\xi\in L^2 (G)$ using the
unimodularity of $g$ \cite[Lemma 5.1]{PS} and that the function
from $L^2 (G) \setminus \{ 0 \}$ to $\mathcal{B} (L^2 (G))$ given
by $\xi\mapsto p_\xi$ is norm continuous. Thus, since for a given
$\xi\in L^2 (G)$ the function from $\Aut (G)$ to $L^2 (G)$ defined
by $\gamma\mapsto\xi\circ\gamma^{-1}$ is continuous by
Proposition~IV.5.2 of \cite{Brac}, we infer in view of our
assumption on $g$ that the set
$$ \Theta = \big\{ \beta_g^n (p_k ) : k=1, \dots , r \text{ and }
n\in\mathbb{Z} \big\} $$
has compact closure in $\mathcal{K} (L^2 (G))$. Now suppose
$I\subseteq\mathbb{Z}$ is a positive density subset and let $\varepsilon > 0$.
Since $\Theta$ has compact closure there is a finite subset
$F\subseteq\Theta$ which is $\delta$-dense in $\Theta$ for $\delta =
\varepsilon\big( 2r \max_{1\leq k \leq r} \| a_k \| \big)^{-1}$.
Using the fact that every finite partition of a positive density subset of
$\mathbb{Z}$ contains at least one member of positive upper density, we can
apply a diagonal argument across
$k$ to find a subset $J\subseteq I$ of positive upper density and
$q_1 , \dots , q_r \in F$ such that
$$ \Big( \max_{1\leq k \leq r} \| a_k \| \Big) \| \beta_g^n (p_k ) - q_k \|
\leq \frac{\varepsilon}{2r} $$
for every $k=1, \dots , r$ and $n\in J$. By Proposition~\ref{P-directsum} the
$\ell_\infty$-direct sum of $r$ copies of $\alpha_g$ has zero CA entropy, and
so by Theorem~\ref{T-zero} there exists a finite subset $E\subseteq J$ and a
norm one element $(c_n )_{n\in E}$ of $\ell_1^E$ such that
$$ \sup_{1\leq k \leq r} \bigg\| \sum_{n\in E} c_n \alpha_g^n (a_k ) \bigg\|
\leq \frac{\varepsilon}{2r} . $$
We then have, for each $k=1, \dots , r$,
\begin{align*}
\lefteqn{\bigg\| \sum_{n\in E} c_n (\alpha_g \otimes\beta_g )^n 
(a_k \otimes p_k ) \bigg\| }\hspace*{26mm} \\
\hspace*{22mm} &\leq \bigg\| \bigg( \sum_{n\in E} c_n \alpha_g^n 
(a_k ) \bigg) \otimes q_k \bigg\| \\
&\hspace*{17mm} \ + \bigg\| \sum_{n\in E} c_n \alpha_g^n (a_k ) \otimes
(\beta_g^n (p_k ) - q_k ) \bigg\| \\
&\leq \bigg\| \sum_{n\in E} c_n \alpha_g^n (a_k ) \bigg\| + \sum_{n\in E}
| c_n | \| a_k \| \| \beta_g^n (p_k ) - q_k \| \\
&\leq \frac{\varepsilon}{2r} + \frac{\varepsilon}{2r} = \frac{\varepsilon}{r}
\end{align*}
and hence
\begin{align*}
\bigg\| \sum_{n\in E} c_n (\alpha_g \otimes\beta_g )^n (x) \bigg\|
&\leq \sum_{k=1}^r \bigg\| \sum_{n\in E} c_n (\alpha_g \otimes\beta_g )^n
(a_k \otimes p_k )\bigg\| \\
&\leq r \frac{\varepsilon}{r} = \varepsilon .
\end{align*}
Since $\varepsilon$ was arbitrary we conclude that $I$
is not an $\ell_1$-isomorphism set for $x$, completing the proof.
\end{proof}

One consequence of Theorem~\ref{T-crossedprod} is the
existence of a simple separable unital nuclear $C^*$-algebra on
which every possible value for CA entropy is realized by an inner
$^*$-automorphism. Indeed we argue that this happens for the
Bunce-Deddens algebra $B$ of type $2^\infty$ (see Example~3.2.11
of \cite{Ror}) as follows. By the classification theorem of
\cite{Ell} $B$ can be expressed as the crossed product associated
to the dyadic odometer, since these two $C^*$-algebras are real
rank zero A$\mathbb{T}$-algebras with the same Elliott invariant
(see Section~3.2 of \cite{Ror}). Now given any $r\in [0,\infty ]$,
by \cite{BH} there is a minimal homeomorphism of the Cantor set
which is strong orbit equivalent to the dyadic odometer and has
topological entropy equal to $r$, and by Theorem~2.1 of \cite{GPS}
the crossed product associated to this homeomorphism is
$^*$-isomorphic to $B$. Applying Theorem~\ref{T-crossedprod}
we thus obtain an inner $^*$-automorphism of $B$ with CA entropy
equal to $r$.

We close this section with a result on reduced free product $^*$-automorphisms.
For notation and terminology see Section~2 of \cite{BDS}.

\begin{proposition}
Let $D$ be a finite-dimensional $C^*$-algebra. Let $(A,E)$ and $(B,F)$ be
$D$-probability spaces with $A$ and $B$ commutative, and suppose that the GNS
representations of $E$ and $F$ are faithful. Let $\alpha$ and $\beta$ be
$^*$-automorphisms of $(A,E)$ and $(B,F)$, respectively, such that
$\alpha |_D = \beta |_D$ and $\hc (\alpha ) = \hc (\beta ) = 0$. Then
$$ \hc (\alpha * \beta ) = 0 . $$
\end{proposition}

\begin{proof}
For a $^*$-automorphism of a unital commutative $C^*$-algebra the CA entropy
and Voiculescu-Brown entropy agree by
Proposition~\ref{P-topol}. We can thus apply Theorem~5.7 of \cite{BDS} to
obtain $\entr (\alpha * \beta ) = 0$, and this implies that
$\hc (\alpha * \beta ) = 0$ by Theorem~\ref{T-zero} and \cite[Thm.\ 3.3]{EID}
(see Proposition~\ref{P-VoiBr}).
\end{proof}

\section{On the prevalence of zero and infinite CA entropy in
$C^*$-algebras}\label{S-opalg}

Let $A$ be a unital $C^*$-algebra $A$. We denote by $\mathcal{U} (A)$ the
unitary group of $A$ and by $\mathcal{U}_0 (A)$ the subgroup of
$\mathcal{U} (A)$ consisting of those unitaries which are homotopic to $1$.
We denote by $\Aut (A)$ the group of $^*$-automorphisms of $A$, by $\Inn (A)$
the subgroup of inner $^*$-automorphisms, and by $\Inn_0 (A)$ the subgroup
of inner $^*$-automorphisms that can be expressed as $\Ad u$ for some
$u\in\mathcal{U}_0 (A)$. Unless stated otherwise, the topology on $\Aut (A)$
will be the point-norm one, i.e., the topology which has as a base
sets of the form $\{ \beta\in\Aut (A) : \| \beta (a) - \alpha (a) \| <
\varepsilon \text{ for all } a\in\Omega \}$ for some $\varepsilon > 0$ and
finite set $\Omega\subseteq\Aut (A)$. In particular, $\overline{\Inn} (A)$ and
$\overline{\Inn_0} (A)$ refer to point-norm closures. For separable $A$
the space $\Aut (A)$ is Polish.

\begin{proposition}\label{P-finitespec}
Let $A$ be a unital $C^*$-algebra and $u\in A$ a unitary with finite
spectrum. Then $\hc (\Ad u) = 0$.
\end{proposition}

\begin{proof}
By the functional calculus there exist pairwise orthogonal projections
$p_1 , \dots , p_n \in A$ with sum $1$
and $\lambda_1 , \dots , \lambda_n \in\mathbb{C}$ of unit modulus such that
$u = \lambda_1 p_1 + \cdots + \lambda_n p_n$. For $1\leq i,j \leq n$
set $A_{ij} = p_i A p_j$. Then $A$ decomposes as the direct sum of the
$A_{ij}$'s, and $\Ad (u)$ acts on $A_{ij}$ by multiplication by
$\lambda_i / \lambda_j$. Thus if $\Omega$ is a finite subset of the union of
the $A_{ij}$'s then $\hc (\Ad u , \Omega ) = 0$. The result now follows
from Proposition~\ref{P-prop}(ii).
\end{proof}

Proposition~\ref{P-finitespec} shows in particular that, within the set
of inner $^*$-automorphisms of a von Neumann algebra, those with zero
CA entropy are norm dense, since any unitary can be approximated
in norm by one with finite spectrum using the Borel functional calculus.

Recall that a $C^*$-algebra is said to be {\it subhomogeneous} if it is
$^*$-isomorphic to a $C^*$-subalgebra of $M_n (C_0 (X))$ for some
$n\in\mathbb{N}$ and locally compact Hausdorff space $X$. By the
structure theory for von Neumann algebras, a von Neumann algebra is
subhomogeneous if and only if it can be written as a finite direct sum of
matrix algebras over commutative von Neumann algebras.

\begin{proposition}
For a von Neumann algebra $M$ the following are equivalent:
\begin{enumerate}
\item $M$ is not subhomogeneous,

\item $\hc (\alpha ) > 0$ for some $\alpha\in\Inn (M)$,

\item The set of $\alpha\in\Inn (M)$ such that $\hc (\alpha ) = \infty$
is norm dense in $\Inn (M)$.
\end{enumerate}
\end{proposition}

\begin{proof}
The implications (3)$\Rightarrow$(1) and (1)$\Leftrightarrow$(2) follow from
Theorem~\ref{T-typeI}. Suppose then that (1) holds. We will obtain (3) upon
showing that, given a unitary $u\in M$
and an $\varepsilon > 0$, there is a unitary $v\in M$
with $\| v-u \| < \varepsilon$ and $\hc (\Ad v) = \infty$.
By the Borel functional calculus we can find a set $\Gamma$ of pairwise
orthogonal projections with sum $1$ in the von Neumann subalgebra of $M$
generated by $u$ with the property that for each $p\in\Gamma$ there is a
$\lambda\in\mathbb{C}$ of unit modulus such that $\| pup - \lambda p \| <
\varepsilon$. By (1) there must be a $q \in\Gamma$ such that $qMq$ is not
subhomogeneous. Then $qMq$ is not a type I $C^*$-algebra and so by
Theorem~\ref{T-typeI} there is a unitary $w\in qMq$ with $\hc(\Ad w) = \infty$.
Perturbing $u$ via the Borel functional calculus if necessary, we may
assume that $quq = \lambda q$ for some $\lambda\in\mathbb{C}$ of unit
modulus. By taking a branch
of the logarithm function we can apply the Borel functional calculus to
obtain a unitary $z\in qMq$ such that
$\| z-1_{qMq} \| < \varepsilon$ and $z^n = w$ for some $n\in\mathbb{N}$. Put
$v=u(z+1-q) \in M$. Then $v$ is a unitary with
$$ \| v-u \| \leq \| z-q \| < \varepsilon . $$
It remains to observe that $\Ad v$ restricts to $\Ad z$ on $qMq$ so that
by monotonicity and Proposition~\ref{P-prop}(iii) we obtain
$$ \hc (\Ad v) \geq \hc (\Ad z) = \frac1n \hc (\Ad z^n ) = \frac1n \hc (\Ad w)
= \infty , $$
completing the proof.
\end{proof}

Since $\Aut (M) = \Inn (M)$ for a type I factor we obtain the following.

\begin{corollary}
Let $M$ be an infinite-dimensional factor of type I. Then the
set of $\alpha\in\Aut (M)$ with $\hc (\alpha) = 0$
(resp.\ $\hc (\alpha) = \infty$) is norm dense in $\Aut (M)$.
\end{corollary}

Let $X$ be a separable Banach space. Choose a dense sequence $\{ x_1 , x_2 ,
\dots \}$ in the unit ball of $X$ and define on $B_1 (X^* )$ the metric
$$ d(\sigma , \omega ) = \sum_{n=1}^\infty 2^{-n}
| \sigma (x_n ) - \omega (x_n ) | , $$
which is compatible with the weak$^*$ topology. On the set of homeomorphisms
of $B_1 (X^* )$ we consider the metric
$$ \rho (T,S) = \sup_{\sigma\in B_1 (X^* )} d(S\sigma , T\sigma ) +
\sup_{\sigma\in B_1 (X^* )} d(S^{-1} \sigma , T^{-1} \sigma ) . $$
Then on $\IA (X)$ the topology arising from $\rho$ (via the identification
an element $\alpha\in\IA (X)$ with the induced homeomorphism $T_\alpha$ of
$B_1 (X^* )$) and the point-norm topology agree, as is readily checked. We can
then apply Lemma 2.4 of \cite{TRP} and Theorem~\ref{T-zero} to obtain the
following lemma.

\begin{lemma}\label{L-Gdelta}
Let $X$ be a separable Banach space. Then the elements $\alpha$ of $\IA (X)$
with $\hc (\alpha ) = 0$ form a $G_\delta$ subset of $\IA (X)$ in
the point-norm topology.
\end{lemma}

\begin{proposition}\label{P-denseGdeltarr0}
Let $A$ be a separable unital $C^*$-algebra with real rank zero. Then
the set of $^*$-automorphisms in $\overline{\Inn_0} (A)$
(resp.\ in $\Inn_0 (A)$) with
zero CA entropy is a dense (resp.\ norm dense) $G_\delta$ subset.
\end{proposition}

\begin{proof}
By \cite{Lin} a unital $C^*$-algebra has real rank zero if and only if
for every unitary $u\in\mathcal{U}_0 (A)$ and $\varepsilon > 0$ there is a
unitary $v\in A$ with finite spectrum such that $\| u-v \| < \varepsilon$.
Thus by Proposition~\ref{P-finitespec} and Lemma~\ref{L-Gdelta} we obtain the
result.
\end{proof}

\begin{proposition}\label{P-denseGdeltaI}
Let $A$ be a separable $C^*$-algebra which is an inductive limit
of type I $C^*$-algebras. Then the set of $^*$-automorphisms in
$\overline{\Inn} (A)$ (resp.\ in $\Inn (A)$) with
zero CA entropy is a dense (resp.\ norm dense) $G_\delta$ subset.
\end{proposition}

\begin{proof}
Apply Theorem~\ref{T-typeI} and Proposition~\ref{P-prop}(ii).
\end{proof}

The following is the noncommutative analogue of Corollary~2.5 in \cite{TRP}.
For the definition of the Jiang-Su $C^*$-algebra $\mathcal{Z}$ see \cite{JS}.

\begin{proposition}\label{P-UHFO2Z}
Let $A$ be a UHF $C^*$-algebra, the Cuntz $C^*$-algebra $\mathcal{O}_2$,
or the Jiang-Su $C^*$-algebra $\mathcal{Z}$.
Then the $^*$-automorphisms in $\Aut (A)$ with zero
CA entropy form a dense $G_\delta$ subset.
\end{proposition}

\begin{proof}
When $A$ is a UHF algebra or $\mathcal{Z}$ then it is an inductive limit of
type I $C^*$-algebras and $\overline{\Inn} (A) = \Aut (A)$ (see \cite{Ror,JS}),
and so the conclusion follows from Proposition~\ref{P-denseGdeltaI}.
In the case of $\mathcal{O}_2$ every unitary is homotopic to $1$ and
every $^*$-automorphism is approximately inner (see \cite{Ror}) whence
$\overline{\Inn_0} (A) = \Aut (A)$, and thus since $\mathcal{O}_2$ has real
rank zero we can apply Proposition~\ref{P-denseGdeltarr0}.
\end{proof}

We will next establish an infinite entropy density result for
$^*$-automorphisms of $C^*$-algebras which are tensorially stable with respect
to the Jiang-Su $C^*$-algebra $\mathcal{Z}$. We say briefly that a
$C^*$-algebra $A$ is {\it $\mathcal{Z}$-stable} if $A\otimes\mathcal{Z}$ is
$^*$-isomorphic to $A$ (note that the $C^*$-tensor product is unique in this
case by the nuclearity of $\mathcal{Z}$). The class of $\mathcal{Z}$-stable
$C^*$-algebras is of importance in the classification theory for nuclear
$C^*$-algebras and includes all Kirchberg algebras and all simple unital
infinite-dimensional AH algebras of bounded dimension \cite[Thm.\ 5]{JS}
\cite[Example 3.4.5]{Ror}  (for obstructions to
$\mathcal{Z}$-stability see \cite{GJS}). Recently A. Toms and W. Winter have
informed us that they have established $\mathcal{Z}$-stability for all
separable approximately divisible $C^*$-algebras.

\begin{lemma}\label{L-aue}
Let $A$ be a $\mathcal{Z}$-stable $C^*$-algebra. For any $\Omega\in
\Fin (A)$ and $\varepsilon > 0$ there exists a $^*$-isomorphism
$\Phi : A\rightarrow A\otimes \mathcal{Z}$ such that
$\| \Phi (x)-x\otimes 1_{\mathcal{Z}} \| < \varepsilon$ for all $x\in \Omega$.
\end{lemma}

\begin{proof}
Since $A$ is $\mathcal{Z}$-stable, we may write $A$ as $A\otimes \mathcal{Z}$,
and by a simple approximation argument we may assume that $\Omega$ is contained
in the algebraic tensor product of $A$ and $\mathcal{Z}$. It suffices then to
consider the case $A = \mathcal{Z}$.
By \cite[Thm.\ 4]{JS} we can find a $^*$-isomorphism $\Psi : \mathcal{Z}
\rightarrow \mathcal{Z}\otimes \mathcal{Z}$. Denote by $\psi$ the embedding
$\mathcal{Z}\rightarrow \mathcal{Z}\otimes \mathcal{Z}$ sending
$a\in \mathcal{Z}$ to $a\otimes 1_{\mathcal{Z}}$. Then $\psi\circ \Psi^{-1}$
is a unital $^*$-endomorphism of $\mathcal{Z}\otimes \mathcal{Z}$ which is
approximately inner by \cite[Thm.\ 3]{JS}. Thus we can find a unitary
$u\in\mathcal{Z}\otimes \mathcal{Z}$ such that
$\| uyu^* - \psi (\Psi^{-1} (y)) \| < \varepsilon$ for all
$y\in \Psi (\Omega )$. Set $\Phi = \Ad u\circ \Psi$.
Then $\| \Phi (x) - x\otimes 1_{\mathcal{Z}} \| < \varepsilon$ for all
$x\in\Omega$.
\end{proof}

\begin{proposition}
Let $A$ be a unital $\mathcal{Z}$-stable $C^*$-algebra. Then the set of
$^*$-automorphisms in $\Aut(A)$ (resp.\ in $\Inn(A)$) with infinite CA
entropy is dense in $\Aut(A)$ (resp.\ in $\Inn(A)$).
\end{proposition}

\begin{proof}
Let $\alpha\in \Aut(A)$, and let $\Omega\in\Fin (A)$ and $\varepsilon > 0$.
By Lemma~\ref{L-aue} we can find a $^*$-isomorphism
$\Phi : A\rightarrow A\otimes \mathcal{Z}$ such that
$\| \Phi (x) - x\otimes 1_{\mathcal{Z}} \| < \varepsilon/2$ for all
$x\in\Omega\cup\alpha (\Omega )$. Since $\mathcal{Z}$ is non-type I, by
Theorem~\ref{T-typeI} it contains a unitary $v$ with $\hc (\Ad v) = \infty$.
Set $\beta = \Phi^{-1} \circ (\alpha\otimes\Ad v)\circ \Phi$.
Note that when $\alpha$ is inner, so is $\beta$. Then
$\hc (\beta ) = \infty$ by monotonicity, and for all $x\in\Omega$ we have
\begin{align*}
\|\beta (x) - \alpha (x)\|
&= \|(\alpha\otimes\Ad v)(\Phi (x))-\Phi (\alpha (x))\| \\
&\leq \|(\alpha\otimes\Ad v)(\Phi (x))-(\alpha\otimes \Ad v)(x\otimes
1_{\mathcal{Z}})\| \\
&\hspace*{20mm} \ + \|\alpha (x)\otimes 1_{\mathcal{Z}} -
\Phi (\alpha (x))\| \\
&< \varepsilon,
\end{align*}
completing the proof.
\end{proof}

\begin{question}
Does there exist a non-type I $C^*$-algebra $A$ such that
the set of $^*$-automorphisms in $\Aut (A)$ with infinite CA entropy is not
dense in $\Aut (A)$?
\end{question}

\begin{question}\label{Q-nuclear}
Let $A$ be a unital nuclear $C^*$-algebra and $u\in A$ a unitary with
$\hc (\Ad u) > 0$. Must the spectrum of $u$ be the entire unit circle?
\end{question}

\noindent Note that the answer to Question~\ref{Q-nuclear} is no if, for
example, $A$ is a von Neumann algebra, for in this case by the Borel functional
calculus there is a unitary $v\in A$ without full spectrum such that $v^2 = u$,
and $\hc (\Ad v) > 0$ by Proposition~\ref{P-prop}(iii).

\section{The shift on $M_p^{\,\otimes\mathbb{Z}}$}\label{S-shift}

In this section we will show that for $p\geq 2$ the tensor product shift on
the UHF $C^*$-algebra $M_p^{\,\otimes\mathbb{Z}}$ (obtained from the shift
$k\mapsto k+1$ on
the index set $\mathbb{Z}$) has infinite CA entropy. The key to obtaining
arbitrarily large lower bounds is the following lemma, which will also
be applied in Sections~\ref{S-ellinfty} and \ref{S-ell1}.

\begin{lemma}\label{L-comb}
Let $X$ be a Banach space and $\alpha\in\IA (X)$. Let $\Omega\subseteq X$ be 
a finite set of unit vectors such that $\| x + y \| < 2$ for all 
$x,y\in\Omega$ with $x\neq y$ and 
$\big\| \sum_{k=0}^{n-1} \alpha^k (x_k ) \big\| = n$
for all $n\in\mathbb{N}$ and $(x_k )_k \in \Omega^n$. Then 
$\hc (\alpha , \Omega ) \geq\log | \Omega |$.
\end{lemma}

\begin{proof}
Let $0 < \delta < 2 - \max \{ \| x + y \| : x,y\in\Omega \text{ and } 
x\neq y \}$ and $n\in\mathbb{N}$. Let $\iota : X\to Y$ be a CA embedding
and let $(\varphi , \psi ,d ) \in\CA (\iota , \Omega\cup\alpha\Omega\cup
\cdots\cup\alpha^{n-1} \Omega , \delta^2 )$. For
every $x = (x_k )_k \in \Omega^n$ we can find a $1\leq b_x \leq d$ such that
$\big| \sum_{k=0}^{n-1} \varphi (\alpha^k (x_k )) (b_x ) \big| \geq
n(1-\delta^2 )$. Take a maximal subset $Q_n \subseteq \Omega^n$ with the
property that $| \{ k : x_k \neq y_k \} | \geq 3n\delta$ for all $x,y\in Q_n$
with $x\neq y$ (in which case $\big\| \sum_{k=0}^{n-1} \alpha^k (x_k + y_k )
\big\| < n(2-3\delta^2 )$). If $x$ and $y$ are any distinct elements
of $Q_n$ with $b_x = b_y$ then the arguments of 
$\sum_{k=0}^{n-1} \varphi (\alpha^k (x_k ))
(b_x )$ and $\sum_{k=0}^{n-1} \varphi (\alpha^k (y_k )) (b_y )$ are
separated by at least $\delta^2$, and so we get the inequality
$d\geq \delta^2 (2\pi )^{-1} | Q_n |$. 

With a view to obtaining a lower bound for $| Q_n |$, for each $x\in Q_n$
define $\Theta_x$ to be the subset of $\Omega^n$ consisting of all
$y$ such that $|\{ k : x_k \neq y_k \} | < 3n\delta$. Setting $m = 
| \Omega |$ and supposing that
$\delta < 1/6$ and $n$ is sufficiently large, we have the estimate
$$ | \Theta_x | = \sum_{k=0}^{\lfloor 3n\delta \rfloor} 
(m-1)^k \binom{n}{k} \leq m^{3n\delta} \binom{n}{3n\delta} . $$
By Stirling's formula we can find an $M>0$ such that for all 
$n\in\mathbb{N}$ and $0 < \delta < 1/6$ we have
$$ \binom{n}{3n\delta} \leq \frac{M}{\sqrt{n\delta}} 
(1-3\delta )^{-n} \left( \frac{1-3\delta }{3\delta} \right)^{3n\delta } . $$
Using the maximality of $Q_n$, it follows in our case that 
\begin{align*}
| Q_n | &\geq \frac{m^n}{\max_{x\in Q_n} | \Theta_x |} \\
&\geq m^{n(1-3\delta )} \frac{\sqrt{n\delta}}{M} (1-3\delta )^n
\left( \frac{3\delta}{1-3\delta } \right)^{3n\delta } . 
\end{align*}
Since the logarithm of this last expression divided by $n$ tends to 
$\log m + C(\delta )$
as $n\to\infty$ where $\lim_{\delta\to 0^+} C(\delta ) = 0$, we obtain the
assertion of the lemma.
\end{proof}

For each $k\in\mathbb{N}$ we denote by $u_k$ and $v_k$
the self-adjoint matrices
$\big[ \begin{smallmatrix} 1 & 0 \\ 0 & -1 \end{smallmatrix} \big]$ and
$\big[ \begin{smallmatrix} 0 & 1 \\ 1 & 0 \end{smallmatrix} \big]$,
respectively, viewed as elements $M_p^{\,\otimes\mathbb{Z}}$
by first identifying $M_2$ with the upper left-hand $2\times 2$ corner of
$M_p$ and then embedding $M_p$ into $M_p^{\,\otimes\mathbb{Z}}$ as the $k$th
tensor product factor. For each
$n\in\mathbb{N}$ we set $Y_n = \spn_{\mathbb{R}} \{ u_k , v_k \}_{k=1}^n
\subseteq (M_p^{\,\otimes\mathbb{Z}})_\sa$ and denote by $W_n$ the
$\ell_1$-direct sum of $n$ copies of $\ell_2^2$ over the real numbers.

\begin{lemma}\label{L-isometric}
For each $n\in\mathbb{N}$ the linear map from $W_n$ to $Y_n$ given by
$$ ((c_k , d_k ))_{k=1}^n \mapsto \sum_{k=1}^n c_k u_k + d_k v_k $$
is an isometric isomorphism.
\end{lemma}

\begin{proof}
Let $c_1 , \dots , c_n , d_1 , \dots , d_n \in\mathbb{R}$. For each
$k=1, \dots , n$ we note that $c_k u_k + d_k v_k$, as an element of $M_2$, has
characteristic polynomial $x^2 - ( c_k^2 + d_k^2 )$ and hence admits
a unit eigenvector whose corresponding vector state $\sigma_k$ satisfies
$$ \sigma_k (c_k u_k + d_k v_k ) = ( c_k^2 + d_k^2 )^{1/2} =
\| c_k u_k + d_k v_k \| . $$
Thus
\begin{align*}
\sum_{k=1}^n ( c_k^2 + d_k^2 )^{1/2} &= 
\sum_{k=1}^n \| c_k u_k + d_k v_k \| \\
&\geq \bigg\| \sum_{k=1}^n c_k u_k + d_k v_k \bigg\| \\
&\geq \bigg| ( \sigma_1 \otimes\cdots\otimes \sigma_n ) \bigg(
\sum_{k=1}^n c_k u_k + d_k v_k \bigg) \bigg| \\
&= \bigg| \sum_{k=1}^n \sigma_k (c_k u_k + d_k v_k ) \bigg| \\
&= \sum_{k=1}^n ( c_k^2 + d_k^2 )^{1/2}
\end{align*}
so that $\big\| \sum_{k=1}^n c_k u_k + d_k v_k \big\| = \sum_{k=1}^n
( c_k^2 + d_k^2 )^{1/2}$, which establishes the proposition.
\end{proof}

\begin{theorem}\label{T-shift}
The tensor product shift $\alpha$ on
$M_p^{\,\otimes\mathbb{Z}}$ has infinite CA entropy.
\end{theorem}

\begin{proof}
By Lemma~\ref{L-isometric} every finite subset of the unit sphere of $Y_1$
satisfies the hypotheses of Lemma~\ref{L-comb}, and so we obtain the result.
\end{proof}

Theorem~\ref{T-shift} prompts the following question. (Compare the discussion
after Theorem~\ref{T-crossedprod}.)

\begin{question}
Is there a $^*$-automorphism of $M_p^{\,\otimes\mathbb{Z}}$ (or any other
simple AF algebra) with finite nonzero CA entropy?
\end{question}

\noindent We also remark that we don't know whether there exists a
$^*$-automorphism of a simple purely infinite $C^*$-algebra (in particular,
the Cuntz algebra $\mathcal{O}_2$) with finite nonzero
CA entropy (cf.\ \cite{BDS}).

\section{Isometric automorphisms of $\ell_\infty$}\label{S-ellinfty}

By the Banach-Stone theorem, for every isometric automorphism $\alpha$ of
$\ell_\infty$ there are a double-sided sequence
$\lambda : \mathbb{Z}\to\mathbb{T}$ and permutation $\sigma$ of $\mathbb{Z}$
such that $\alpha (x)(s) = \lambda (s) x(\sigma (s))$ for all
$x\in\ell_\infty$ and $s\in\mathbb{Z}$ (cf.\ \cite[Prop.\ 2.f.14]{CBS}).

\begin{proposition}\label{P-ellinftyIA}
For every isometric automorphism $\alpha$ of $\ell_\infty$ we have either
$\hc (\alpha ) = 0$ or $\hc (\alpha ) = \infty$ depending on whether or
not there is a finite bound on the cardinality of the orbits of
the associated permutation of $\mathbb{Z}$.
\end{proposition}

\begin{proof}
Let $\alpha$ be an isometric automorphism of $\ell_\infty$ with associated
$\mathbb{T}$-valued sequence $\lambda$ and permutation
$\sigma$ of $\mathbb{Z}$, as above. Suppose first that there is a
$d\in\mathbb{N}$ such that every orbit of $\sigma$ is of cardinality at
most $d$. To show that $\hc (\alpha ) = 0$ we may assume that $\sigma$
is the identity by replacing $\alpha$ with $\alpha^{d!}$ and applying
Proposition~\ref{P-prop}(iii). We may then view $\alpha$ as the isometric
automorphism of $C(\beta\mathbb{Z} ) \cong\ell_\infty$ given by
multiplication by the continuous $\mathbb{T}$-valued extension of $\lambda$ to
the Stone-\v{C}ech compactification $\beta\mathbb{Z}$, so that we can appeal
to Proposition~\ref{P-cms} to conclude that $\hc (\alpha ) = 0$,
as desired.

Suppose now that there is no finite bound on the
cardinality of the orbits of $\sigma$. Let $m\in\mathbb{N}$. 
By relabelling the coordinates of
$\ell_\infty$ if necessary we can associate to each $n\in\mathbb{N}$ and
$f\in\{ 1, \dots , m \}^{\{ 0, \dots , n-1 \}}$
an interval $I_f$ in $\mathbb{Z}$ of the form $[s_f , s_f + n-1 ]$ such that
$\sigma (t) = t+1$ for all $t\in I_f$, and we can arrange for the intervals
$I_f$ to be pairwise disjoint. Now
let $1\leq j\leq m$ and define the unit vector $x_j \in\ell_\infty$
as follows. For every $n\in\mathbb{N}$,
$f\in \{ 1, \dots , m \}^{\{ 0, \dots , n-1 \}}$, and $k=0, \dots ,n-1$ we set
$$ x_j (s_f + k) = \left\{ \begin{array}{l@{\hspace*{8mm}}l} 
\overline{\lambda (s_f ) \lambda (s_f + 1)\cdots\lambda (s_f + k - 1)} & 
\text{if } f(k) = j, \\ 
0 & \text{otherwise}.
\end{array} \right. $$
At all remaining coordinates we take $x_j$ to be zero. Then 
the set $\Omega = \{ x_1 , \dots , x_m \}$ satisfies the hypotheses of 
Lemma~\ref{L-comb}, from which we obtain $\hc (\alpha , \Omega ) \geq
\log m$. Since $m$ is arbitrary we conclude that $\hc (\alpha ) = \infty$,
completing the proof.
\end{proof}

\begin{remark}
From Propositions~\ref{P-topol} and \ref{P-ellinftyIA} we see that the
Stone-\v{C}ech compactification $\beta\mathbb{Z}$ (i.e., the
spectrum of $\ell_\infty$) provides an example of a compact Hausdorff space
which admits a homeomorphism with infinite topological entropy
but no homeomorphism with finite nonzero topological entropy.
\end{remark}

\section{Isometric automorphisms of $\ell_1$}\label{S-ell1}

Given an isometric automorphism $\alpha$ of
$\ell_1$ there are a double-sided sequence
$\lambda : \mathbb{Z}\to\mathbb{T}$ and a permutation $\sigma$ of $\mathbb{Z}$
such that $\alpha (x)(s) = \lambda (s) x(\sigma (s))$ for all
$x\in\ell_1$ and $s\in\mathbb{Z}$ \cite[Prop.\ 2.f.14]{CBS}.

\begin{proposition}\label{P-ell1}
For every isometric automorphism $\alpha$ of $\ell_1$ we have either
$\hc (\alpha ) = \infty$ or $\hc (\alpha ) = 0$ depending on whether or
not there is an infinite orbit in
the associated permutation of $\mathbb{Z}$.
\end{proposition}

\begin{proof}
Suppose first that there is an infinite orbit in the permutation of
$\mathbb{Z}$ associated to $\alpha$. Let $k$ be a integer contained in this
infinite orbit and let $e_k$ be the $k$th standard unit basis vector
in $\ell_1$. Then for any finite subset $\Omega$ of the unit sphere of
$\spn \{ e_k \}$ we have $\hc (\alpha , \Omega ) \geq\log | \Omega |$
by Lemma~\ref{L-comb}, whence $\hc (\alpha ) = \infty$.

If on the other hand there is no 
infinite orbit in the associated permutation of
$\mathbb{Z}$, then $\ell_1$ contains a dense union of finite-dimensional
$\alpha$-invariant subspaces, so that $\hc (\alpha ) = 0$ by
Proposition~\ref{P-prop}(ii).
\end{proof}

\begin{remark}\label{R-sphere}
The infinite value of CA entropy for the shift on $\ell_1$ is a reflection of
the fact that the unit sphere of the scalar field $\mathbb{C}$ is infinite.
If we work instead over $\mathbb{R}$, whose unit sphere has two elements, then
the CA entropy of the shift $\alpha'$ on $(\ell_1)_{\mathbb{R}}$ is $\log 2$.
Indeed $(\ell^n_1 )_{\mathbb{R}}$ embeds isometrically in
$(\ell^{2^n}_\infty )_{\mathbb{R}}$ so that $\hc (\alpha' ) \leq\log 2$,
while from Lemma~\ref{L-comb} we get the lower bound
$\hc (\alpha' ) \geq\log 2$ by taking $\Omega = \{ e , -e \}$, where 
$e$ is any standard unit basis vector in $(\ell_1 )_{\mathbb{R}}$.
\end{remark}

\begin{remark}\label{R-conj}
Using Proposition~\ref{P-ell1} we can show that infinite CA entropy is
not an isomorphic conjugacy invariant, in contrast to zero CA entropy.
Consider the shift $T$ on
$\{ 1,-1 \}^{\mathbb{Z}}$. Define the self-adjoint unitary function
$f\in C(\{ 1,-1 \}^{\mathbb{Z}} )$ by
$$ f( ( a_k )_k ) = a_0 $$
for all $( a_k )_k \in\{ 1,-1 \}^{\mathbb{Z}}$. Then
$\{ f\circ T^k \}_{k\in\mathbb{Z}}$ is isometrically equivalent to the standard
basis of $\ell_1$ over $\mathbb{R}$, and hence equivalent to the standard
basis of $\ell_1$ over $\mathbb{C}$. Denoting by $\beta_T$ the restriction of
the induced $^*$-automorphism of $C(\{ 1,-1 \}^{\mathbb{Z}} )$ to
$\overline{\spn} \{ f\circ T^k \}_{k\in\mathbb{Z}}$, we obtain an isometric
automorphism which is isomorphically conjugate to the shift on $\ell_1$.
While the latter has infinite CA entropy, however, we have
$\hc (\beta_T ) \leq \log 2$ by monotonicity and Proposition~\ref{P-topol}.
\end{remark}

\begin{example}\label{E-free}
Let $\sigma$ be a permutation of $\mathbb{Z}$. Let $\sigma_*$ be the
corresponding free permutation $^*$-automorphism of the full free group
$C^*$-algebra $C^* (F_{\mathbb{Z}})$ sending $u_j$ to $u_{\sigma(j)}$, where
$\{u_j\}_{j\in\mathbb{Z}}$ is the set of canonical unitary generators of
$C^* (F_{\mathbb{Z}})$. Then $\hc (\sigma_* ) = \infty$ if and only
if $\sigma$ has an infinite orbit. As pointed out in \cite[Sect.\ 8]{Pis},
$\{u_j\}_{j\in\mathbb{Z}}$ is isometrically equivalent to the standard basis
of $\ell_1$, and so the ``if'' part holds by Proposition~\ref{P-ell1} and
monotonicity. For the ``only if'' part, supposing that
$\sigma$ has no infinite orbits we can take $\Delta$ in Theorem~\ref{T-zero}
to be the set of unitaries in $C^*(F_{\mathbb{Z}})$ corresponding to elements
of $F_{\mathbb{Z}}$ to get $\htopol (T_{\sigma_*}) = 0$, as desired.
We remark that the corresponding $^*$-automorphisms of the reduced free group
$C^*$-algebra $C_\red^* (F_{\mathbb{Z}})$ all have zero CA entropy by
\cite{BC,Dyk} and Proposition~\ref{P-VoiBr}.
\end{example}

\end{document}